\documentclass[english]{IEEEtran}
\usepackage[T1]{fontenc}
\usepackage[utf8]{inputenc}
\usepackage{array}
\usepackage{float}
\usepackage{mathtools}
\usepackage{multirow}
\usepackage{amsmath}
\usepackage{amsthm}
\usepackage{amssymb}
\usepackage{graphicx}

\makeatletter

\providecommand{\tabularnewline}{\\}
\floatstyle{ruled}
\newfloat{algorithm}{tbp}{loa}
\providecommand{\algorithmname}{Algorithm}
\floatname{algorithm}{\protect\algorithmname}

\theoremstyle{plain}
\newtheorem{thm}{\protect\theoremname}
\theoremstyle{plain}
\newtheorem{prop}[thm]{\protect\propositionname}

\usepackage{fontawesome}
\usepackage{color}
\usepackage{cite}
\usepackage{amsthm}
\theoremstyle{plain}
\newtheorem{asm}{Assumption}
\newtheorem*{asm*}{Assumption}
\usepackage{xcolor}
\usepackage{multirow}

\makeatother

\usepackage{babel}
\providecommand{\propositionname}{Proposition}
\providecommand{\theoremname}{Theorem}

\begin{document}
\title{Power System State Estimation by Phase Synchronization and Eigenvectors}
\author{Iven Guzel, \emph{Student Member, IEEE}, Richard Y. Zhang,\emph{ Member,
IEEE}\thanks{The authors are with the Department of Electrical and Computer Engineering,
University of Illinois Urbana Champaign, Urbana, IL 61801 USA (e-mails:
iguzel2@illinois.edu, ryz@illinois.edu).}}

\maketitle
\global\long\def\R{\mathbb{R}}%
\global\long\def\C{\mathbb{\mathbb{C}}}%
\global\long\def\Aop{\mathcal{A}}%
\global\long\def\Norm{\mathrm{Norm}}%
\global\long\def\gnd{\mathrm{gnd}}%
\global\long\def\opt{\mathrm{opt}}%
\global\long\def\bus{\mathrm{bus}}%
\global\long\def\From{\mathrm{f}}%
\global\long\def\To{\mathrm{t}}%
\global\long\def\lb{\mathrm{lb}}%
\global\long\def\one{\mathbf{\mathbf{1}}}%
\global\long\def\half{{\textstyle \frac{1}{2}}}%
\global\long\def\e{\mathbf{e}}%
\global\long\def\i{\mathrm{j}}%
\global\long\def\inner#1#2{\left\langle #1,#2\right\rangle }%
\global\long\def\o#1{\overline{#1}}%
\global\long\def\Re{\operatorname{Re}}%

\global\long\def\Im{\operatorname{Im}}%
\global\long\def\tr{\operatorname{tr}}%
\global\long\def\rank{\operatorname{rank}}%
\global\long\def\diag{\operatorname{diag}}%
\global\long\def\proj{\operatorname{proj}}%
\global\long\def\conj{\operatorname{conj}}%
\global\long\def\H{\mathbb{H}}%
\global\long\def\pmu{\mathrm{pmu}}%

\begin{abstract}
To estimate accurate voltage phasors from inaccurate voltage magnitude
and complex power measurements, the standard approach is to iteratively
refine a good initial guess using the Gauss--Newton method. But the
nonconvexity of the estimation makes the Gauss--Newton method sensitive
to its initial guess, so human intervention is needed to detect convergence
to plausible but ultimately spurious estimates. This paper makes a
novel connection between the angle estimation subproblem and phase
synchronization to yield two key benefits: (1) an exceptionally high
quality initial guess over the angles, known as a \emph{spectral initialization};
(2) a correctness guarantee for the estimated angles, known as a \emph{global
optimality certificate}. These are formulated as sparse eigenvalue-eigenvector
problems, which we efficiently compute in time comparable to a few
Gauss-Newton iterations. Our experiments on the complete set of Polish,
PEGASE, and RTE models show, where voltage magnitudes are already
reasonably accurate, that spectral initialization provides an almost-perfect
single-shot estimation of $n$ angles from $2n$ moderately noisy
bus power measurements (i.e. $n$ pairs of PQ measurements), whose
correctness becomes guaranteed after a single Gauss--Newton iteration.
For less accurate voltage magnitudes, the performance of the method
degrades gracefully; even with moderate voltage magnitude errors,
the estimated voltage angles remain surprisingly accurate.
\end{abstract}

\section{Introduction }\label{sec:intro}

Power system state estimation (PSSE) is a critical component of grid
operations, vital for achieving and maintaining situational awareness
in an increasingly variable, uncertain, and automated grid. The basic
function of PSSE is to estimate a complete set of complex voltage
phasors (i.e. magnitudes and angles at each bus) from a potentially
incomplete set of noisy voltage magnitude and complex power measurements.
To achieve this, the standard approach is to maximize a certain likelihood
function, by iteratively refining an initial guess using the Gauss--Newton
method, in order to compute the maximum likelihood estimator (MLE)~\cite{schweppe1970power}.

Unfortunately, the likelihood function for PSSE is highly nonlinear
and nonconvex, and this makes iterative refinement very sensitive
to the quality of the initial guess. In theory, if the initial guess
is sufficiently close to the MLE, then rapid and numerically robust
convergence is guaranteed (see e.g. \cite[Theorem 10.1]{nocedal2006numerical}
or \cite[Theorem~10.2.1]{dennis1996numerical}), so long as the system
is observable and the measurements are reasonably accurate. But in
practice, even when using high-quality initializations based on prior
estimates, convergence can still be slow and numerically unstable.
This fundamental gap arises because even a very good initial guess
from a practical perspective is rarely accurate enough for theoretical
guarantees to be rigorously applicable. The resulting convergence
issues tend to become more severe as the system is placed under stress,
and when it is subject to sudden changes. In fact, convergence can
fail altogether during unusual or emergency situations, precisely
when accurate state estimation is needed the most~\cite{zhao2018statistical}.

\begin{figure}[t]
\includegraphics[width=1\columnwidth]{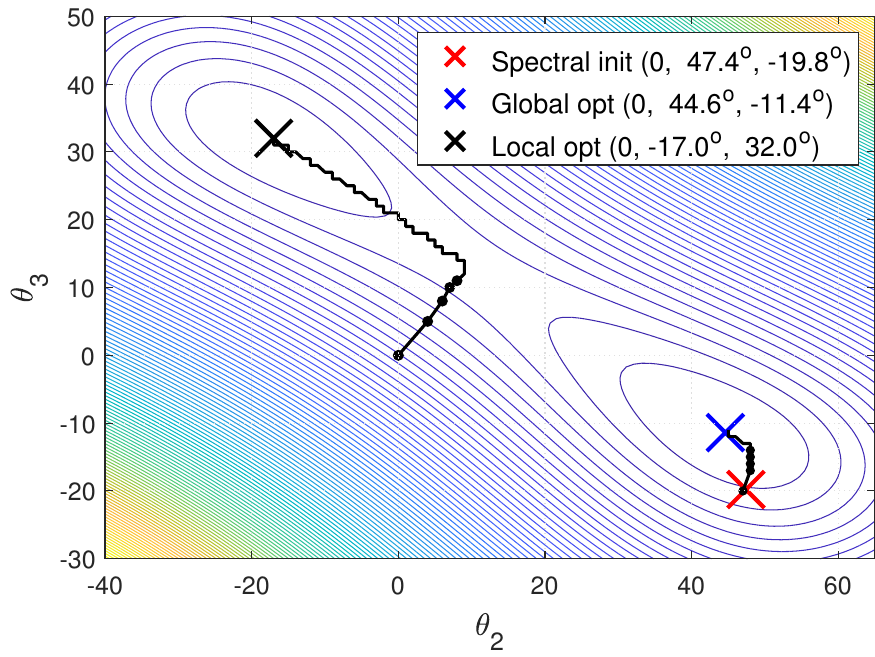}

\vspace{-0.3cm}

\caption{\textbf{Gauss--Newton can produce misleading spurious estimates over
the voltage angles, even when the voltage magnitudes are perfectly
measured}. For a three-bus system with perfect voltage magnitude measurements
(see Appendix A), the usual ``cold start'' initialization experiences
convergence issues, resulting in slow convergence over many iterations
to a spurious estimate ({\textcolor{black}{\faTimes}}). In contrast,
our spectral initialization technique gives a one-shot estimate ({\textcolor{red}{\faTimes}})
close to the maximum likelihood estimator ({\textcolor{blue}{\faTimes}})
in a region where the convergence behavior is vastly improved. Dots
($\bullet$) show the first 5 Gauss--Newton iterations.}

\label{fig:3bus}\vspace{-0.4cm}
\end{figure}

Even when iterative refinement does converge, there is usually no
way to tell whether it actually converged to the MLE. The best we
can say is that the estimator locally maximizes its likelihood, but
local optimality does not imply global optimality in a nonconvex setting.
The danger is convergence to a \emph{spurious} local maximum that
corresponds to a plausible but ultimately fictitious estimate of the
true system state~\cite{zhang2019spurious}. Counterintuitively,
such spurious estimates can be unrelated to the classical phenomenon
of multiple power flow solutions~\cite[Section III-C]{zhang2019spurious}.
As shown in Figure~\ref{fig:3bus}, spurious estimates can arise
even in systems with perfect voltage magnitude measurements, which
by construction cannot exhibit distinct high-voltage and low-voltage
solutions. Currently, this vulnerability to spurious estimates is
addressed by human verification, but it would be desirable to have
an algorithmic means of doing so.

This paper introduces a novel connection between the angle estimation
subproblem of PSSE and a fundamental problem in statistical inference
and signal processing called \emph{phase synchronization}~\cite{singer2011angular,boumal2016nonconvex,bandeira2017tightness}.
An immediate benefit is a high-quality initial guess over the angles
known as a \emph{spectral initialization}. Our experiments on the
complete set of Polish, PEGASE, and RTE models from the IEEE PES PGLib-OPF
dataset~\cite{babaeinejadsarookolaee2019power} find that spectral
initialization is \emph{exceptionally accurate} for observable PSSE
instances with reasonably accurate power and magnitude measurements,
to the extent that it can be used without further refinement in most
practical situations. In the limit, given near-perfect measurements,
spectral initialization outputs the near-exact MLE in one shot. As
such, spectral initialization can be used to compute an initial guess
close enough to the MLE such that rapid and robust convergence can
be mathematically guaranteed. This excellent accuracy degrades gracefully
as the measurements become less accurate; even with moderately inaccurate
measurements, spectral initialization can still remain surprisingly
accurate. 

A second benefit is the ability to \emph{certify the global optimality}
of the estimated angles with respect to a given set of fixed voltage
magnitudes. This can be used to verify the correctness of the estimated
angles, conditioned on the belief that the estimated magnitudes are
accurate. Up to moderate levels of noise, we observe that the certification
comes with \emph{zero duality gap}, so if a given set of angles globally
maximize the likelihood function for fixed magnitudes, then this fact
can be rigorously justified with no conservatism. In practice, we
find that spectral initialization followed by just a single Gauss--Newton
iteration is often enough to estimate angles with certified global
optimality.  Although our certificate comes with a conditional clause
with respect to the voltage magnitudes, it is the first nonvacuous
certificate for PSSE under any restriction, to the best of our knowledge. 

The time it takes to compute the spectral initialization and the global
optimality certification poses a potentially serious concern. Both
require computing eigenvalue-eigenvector pairs for an $n\times n$
complex Hermitian matrix, where the number of buses $n$ could be
as large as tens of thousands. Our critical insight is that the \emph{chordal
sparsity} (and not just mere sparsity) of power systems~\cite{zhang2023parameterized}
allows us to solve these eigenproblems using the inverse iteration~\cite{ipsen1997computing},
which is substantially more convergent, more robust, and more accurate
than the (generalized) power iterations typically used for phase synchronization~\cite{boumal2016nonconvex}.
Using the inverse iteration, we consistently and reliably compute
the spectral initialization and the global optimality certification
in time comparable to a few Gauss--Newton iterations.

Combined, our approach helps to address the sensitivity of iterative
refinement to the quality of the initial guess, for cost equivalent
to a few extra iterations. Its main limitation, however, is that it
sidesteps the voltage magnitudes; it does not provide an initial guess
for the magnitudes, and is only able to certify the global optimality
of angles subject to fixed magnitudes. Therefore, our proposed approach
is most effective when a good initial guess is already available for
the voltage magnitudes. This is the case if most bus voltage magnitudes
are directly measured by SCADA telemetry, and that any missing voltage
magnitudes can be assumed to lie within a known range. We argue that
this is a reasonable compromise, given that the original impetus and
the core difficulty of PSSE is to statistically infer voltage angles
because they cannot be directly measured~\cite{schweppe1970power}.

Our approach is most closely related to prior work using semidefinite
programming (SDP) relaxations~\cite{zhu2014power,zhang2017conic}
to find good initializations and certify global optimality. However,
traditional SDP relaxations for the PSSE problem tend to be very loose
in the presence of noisy measurements. Consequently, the initializations
they produce can be far from the MLE, and the corresponding lower
bounds are usually vacuous. In contrast, our method can be understood
as solving an underlying SDP relaxation for the phase synchronization
problem, which fixes the voltage magnitudes and focuses solely on
estimating the angles. It is precisely by fixing the magnitudes that
the SDP relaxation becomes tight for estimating angles from moderately
noisy power measurements. To the best of our knowledge, this is the
first tight SDP relaxation for PSSE with noisy measurements, yielding
the first nonvacuous global optimality certification under any restrictions. 

Our focus on estimating provably correct angles from power measurements
gives our approach a superficial similarity with DC estimation~\cite{abur2004power}.
The critical difference is that our proposed approach maintains the
trigonometric nonlinearity between power and angles. It achieves provable
correctness on the nonconvex angle estimation problem subject to fixed
magnitudes with \emph{no approximation}. The ability to handle this
nonlinearity is critical in situations with large angular differences,
when the system is most at risk of collapse.

Besides estimating states from reasonably accurate measurements, another
concern is the potential existence of \emph{bad data}~\cite{merrill1971bad},
i.e. outlier measurements that are grossly inaccurate. The problem
with bad data is that they significantly bias the MLE, and also considerably
worsen convergence. In the existing literature, bad data is typically
handled in one of two approaches: (i) by identifying and removing
the bad data through the classical MLE; (ii) by adopting a robust
MLE capable of rejecting bad data. Consistent with prior work on phase
synchronization~\cite{singer2011angular}, our method is able to
compute the classical MLE even in the presence of bad data. As such,
it can be extended to handle bad data through the first approach of
identifying and removing bad data, though we leave this extension
as future work. 

Finally, compared with prior work on phase synchronization, our work
has some important differences. The classical formulation estimates
$n$ unknown phases from up to $n^{2}$ noisy measurements of pairwise
phase differences $\theta_{i}-\theta_{j}$~\cite{singer2011angular,boumal2016nonconvex},
but our formulation estimates them from a smaller number of \emph{weighted
sums} of pairwise phase differences, e.g. $n$ pairs of complex bus
power measurements. Therefore, the wealth of theoretical guarantees
for phase synchronization do not extent to our novel formulation of
PSSE. It remains important future work to provide a theoretical foundation
for the empirical results reported in this paper.

\textbf{Notation. } The sets $\R^{n}$, $\C^{n}$, and $\mathbb{H}^{n}$
are the real and complex length-$n$ vectors, and the complex Hermitian
$n\times n$ matrices. Subscripts denote indexing, the superscripts
$T$ and $*$ refer to the transpose and Hermitian transpose respectively.
Imaginary unit is $\i\equiv\sqrt{-1}$.  The elementwise product
is denoted $\odot$. The $i$-th standard basis vector is $\e_{i}$.
 We use $X\succeq0$ (resp. $X\succ0$) to express $X\in\mathbb{H}^{n}$
is positive semidefinite (resp. positive definite). We use $\lambda_{\min}(X)$
to refer to minimum eigenvalue of $X$. The diagonal matrix $X$ having
$X_{ii}=x_{i},\ i=1,\dots,n$ is denoted by $\diag(x)$.

\section{Background}\label{sec:backgrnd}

\subsection{System Modeling and Quadratic Measurements}

Consider a system with $n$ buses (vertices) and $l$ branches (edges)
having the bus voltage and bus current injection phasors $v\in\mathbb{C}^{n}$
and $i\in\mathbb{C}^{n}$. We let $i_{\From}\in\mathbb{C}^{l}$ and
$i_{\To}\in\mathbb{C}^{l}$ be the vectors of branch currents at the
\textit{from} and \textit{to} ends of all branches, noting that it
is possible to have multiple branches between two buses. We model
all transmission lines and transformers by a common $\pi$-equivalent
line model, with series admittance $y=(r+\i x)^{-1}$ and total charging
susceptance $b$. Based on this assumption, the network is electrically
linear, and Ohm's Law relates $v$ to bus current injection $i$ and
branch current injections $i_{\From}$ and $i_{\To}$ as follows 
\begin{equation}
i=Y_{\bus}v,\quad i_{\From}=Y_{\From}v,\quad i_{\To}=Y_{\To}v,\label{eq:i_is_Yv}
\end{equation}
where $Y_{\bus}\in\mathbb{C}^{n\times n}$, $Y_{\From}\in\mathbb{C}^{l\times n}$,
and $Y_{\To}\in\mathbb{C}^{l\times n}$ are the bus, \textit{from}
and \textit{to} branch admittance matrices, respectively; see~\cite{zimmerman2010matpower}
for the precise formulas for $Y_{\bus},Y_{\From},Y_{\To}$. We use
per unit notation in all measurements and parameters.

The \emph{squared} voltage magnitude at the $k$-th bus is a quadratic
function of complex bus voltage vector $v\in\C^{n}$, as in 
\begin{equation}
|v_{k}|^{2}=(\e_{k}^{T}v)^{*}(\e_{k}^{T}v)=v^{*}\e_{k}\e_{k}^{T}v.\label{eq:quadratic_v}
\end{equation}
Assuming an electrically linear network (\ref{eq:i_is_Yv}), the complex
bus injection at the $k$-th bus can also be written as a quadratic
function of $v$ 
\begin{equation}
s_{k}=v_{k}i_{k}^{*}=(\e_{k}^{T}i)^{*}(\e_{k}^{T}v)=v^{*}Y_{\bus}^{*}\e_{k}\e_{k}^{T}v.\label{eq:S_k}
\end{equation}
Similarly, the complex power flows on the $k$-th branch spanning
from bus $\alpha(k)$ to bus $\beta(k)$ are also quadratic 
\begin{equation}
\begin{aligned}s_{k,\From} & =v_{\alpha(k)}i_{k,\From}^{*}=(\e_{k}^{T}i_{\From})^{*}(\e_{\alpha(k)}^{T}v)=v^{*}Y_{\From}^{*}\e_{k}\e_{\alpha(k)}^{T}v,\\
s_{k,\To} & =v_{\beta(k)}i_{k,\To}^{*}=(\e_{k}^{T}i_{\To})^{*}(\e_{\beta(k)}^{T}v)=v^{*}Y_{\To}^{*}\e_{k}\e_{\beta(k)}^{T}v.
\end{aligned}
\label{eq:S_flow}
\end{equation}

\subsection{PSSE via Gauss--Newton Method}

The standard formulation of the PSSE problem seeks to estimate the
unknown ground truth vector $v_{\gnd}\in\C^{n}$ given $m$ quadratic
measurements of the following form 
\begin{equation}
b_{i}=v_{\gnd}^{*}A_{i}v_{\gnd}+\epsilon_{i},\qquad\text{where }\epsilon_{i}\overset{\text{i.i.d.}}{\sim}\mathcal{N}(0,\sigma_{i}^{2}).\label{eq:b_i}
\end{equation}
Critically, each $i$-th measurement error $\epsilon_{i}$ is assumed
to be distributed independently and identically according to the zero-mean
Gaussian with standard deviation $\sigma_{i}$. The matrices $A_{i}\in\mathbb{H}^{n}$
relating to voltage magnitude squared, the real and reactive power
injections at the $k$-th bus can be computed from (\ref{eq:quadratic_v})
and (\ref{eq:S_k}). Similarly, the matrices $A_{i}\in\mathbb{H}^{n}$
associated with the real and reactive powers injected into the \textit{from}
and \textit{to} nodes of the $k$-th branch are computed using~(\ref{eq:S_flow}).

The usual way to estimate $v_{\gnd}$ from the noisy measurements
in (\ref{eq:b_i}) is to compute the maximum likelihood estimator
(MLE), which is defined as the global minimizer to the following nonconvex
nonlinear least-squares cost function (corresponding to the negative
log likelihood) 
\begin{equation}
\hat{v}=\arg\min_{v\in\C^{n}}\sum_{i=1}^{m}w_{i}|v^{*}A_{i}v-b_{i}|^{2},\label{eq:psse}
\end{equation}
with weights $w_{i}=\sigma_{i}^{-2}$. In practice, (\ref{eq:psse})
is most commonly solved using the Gauss--Newton method with a polar
parameterization $v_{i}=u_{i}\exp(\i\theta_{i})$. We first rewrite
(\ref{eq:psse}) as the following 
\[
(\hat{\theta},\hat{u})=\underset{(\theta,u)\in\R^{2n-1}}{\arg\min}\|r(\theta,u)\|^{2}
\]
in which one of the phasor angles is held fixed as reference. Beginning
with an initial estimate $v\in\C^{n}$, each iteration computes a
refined estimate $v_{+}\in\C^{n}$ by solving a linear least-squares
problem
\begin{equation}
(\theta_{+},u_{+})=\underset{(\theta_{+},u_{+})\in\R^{2n-1}}{\arg\min}\left\Vert r(\theta,u)+[J(\theta,u)]\begin{bmatrix}\theta_{+}-\theta\\
u_{+}-u
\end{bmatrix}\right\Vert ^{2}\label{eq:gn_iterate}
\end{equation}
where $J(\theta,u)$ is the Jacobian associated with $r(\theta,u)$.
If the initial estimate $v$ is sufficiently close to the global minimum
$\hat{v}$, then the updated estimate $v_{+}$ will converge linearly
towards $\hat{v}$. Below, $\|\cdot\|$ denotes the maximum singular
value of a matrix.
\begin{thm}[Local linear convergence]
\label{thm:lin}Define $\hat{v}$ as the global minimum to (\ref{eq:psse}),
denote $\hat{J}\equiv J(\hat{\theta},\hat{u})$ as the Jacobian evaluated
at $\hat{v}$. If $\sigma<\lambda/2$, where
\[
\lambda\equiv\lambda_{\min}(\hat{J}^{T}\hat{J}),\quad\sigma\equiv2\sigma_{\max}\left(\sum_{i=1}^{m}w_{i}[\hat{v}^{*}A_{i}\hat{v}-b_{i}]A_{i}\right),
\]
then there exists a ball of radius $\rho>0$ about $\hat{v}$ such
that 
\begin{align*}
\|v-\hat{v}\|\le\rho\implies & \|v_{+}-\hat{v}\|\le\frac{2\sigma}{\lambda}\|v-\hat{v}\|<\|v-\hat{v}\|,\\
 & \lambda_{\min}(J(u,\theta)^{T}J(u,\theta))\ge\lambda/2,
\end{align*}
where $v$ and $v_{+}$ are the initial and refined estimates in (\ref{eq:gn_iterate})
respectively.
\end{thm}
\begin{IEEEproof}
This is a direct application of \cite[Theorem~10.2.1]{dennis1996numerical}
with constant $c=2$. 
\end{IEEEproof}
Note that $\lambda$ is a measure of \emph{observability} for the
PSSE problem (see~\cite{fetzer1975observability}), whereas $\sigma$
is proportional to the level of noise in the measurements (\ref{eq:b_i}).
Accordingly, Theorem~\ref{thm:lin} guarantees rapid local convergence
assuming that the system is observable, and that the level of noise
is below a threshold dictated by the observability of the problem.
Moreover, this convergence is guaranteed to be numerically robust,
because the Jacobian matrix $J(\theta,u)$ cannot become ill-conditioned. 

Unfortunately, nothing is guaranteed if Gauss--Newton is initialized
outside of this local neighborhood. In practice, numerical issues
can often cause slow convergence or divergence, even when a very good
initial estimate $v$ is used. While this poor convergence can be
addressed to an extent by adopting backtracking line search, or a
trust-region strategy like the Levenberg--Marquardt algorithm, the
best that can be achieved is a slight enlargement of the local neighborhood
where rapid and numerically robust convergence is guaranteed.

\section{PSSE as Phase Synchronization Problem}\label{sec:reformulation}

In order to make spectral initialization and optimality certification
available for PSSE, we begin by rewriting (\ref{eq:psse}) as a two-stage
problem, with an outer estimation problem over the voltage magnitudes
\begin{subequations}
\label{eq:psse_subprob} 
\begin{gather}
\min_{u\in\R_{+}^{n}}\opt(u),\label{eq:psse_subprob_u}
\end{gather}
and an inner estimation problem over the voltage angles with fixed
voltage magnitudes 
\begin{gather}
\begin{aligned}\opt(u)=\min_{x\in\C^{n}} & \sum_{i=1}^{m}w_{i}|(u\odot x)^{*}A_{i}(u\odot x)-b_{i}|^{2}\\
\text{s.t. } & |x_{i}|=1\text{ for all }i\in\{1,\dots,n\}.
\end{aligned}
\label{eq:psse_subprob_x}
\end{gather}
\end{subequations}
Recall that $\odot$ denotes the elementwise Hadamard product, so
$v=u\odot x$ means $v_{i}=u_{i}x_{i}$ for all $i$. Obviously, the
solution $u$ to the outer subproblem (\ref{eq:psse_subprob_u}) is
exactly $\hat{u}=|\hat{v}|$, and the solution $\hat{x}$ to the inner
subproblem (\ref{eq:psse_subprob_x}) in turn recovers $\hat{v}=\hat{u}\odot\hat{x}$.

For practical reasons, the outer estimation problem (\ref{eq:psse_subprob_u})
over the voltage magnitudes $u$ is generally easier to solve. The
grid operator is required to keep each voltage magnitude $u_{i}$
within a narrow band, such that the ``cold start'' estimation of
$u=\one$ per unit can be expected to have a per-bus deviation of
no more than 0.1 per unit across almost all buses, except perhaps
a very small number of buses at fault. Moreover, in modern power systems,
the vast majority of bus voltage magnitudes are explicitly measured,
so that highly accurate values of $u$ are often explicitly known.

In contrast, it is the inner estimation problem (\ref{eq:psse_subprob_x})
over the voltage angles $\theta$ for $\exp(\i\theta_{i})=x_{i}$
that renders the overall PSSE problem difficult to solve. In this
paper, our crucial insight is to \emph{exactly} reformulate the angle
subproblem (\ref{eq:psse_subprob_x}) into an instance of the phase
synchronization problem, which reads 
\begin{align}
\opt(u)=\min_{x\in\C^{n}}\quad & x^{*}H(u)x\label{eq:psse_ps_u}\\
\text{s.t. }\quad & |x_{i}|=1\text{ for all }i\in\{1,\dots,n\}\nonumber 
\end{align}
where $H$ is a matrix-valued function that we will specify throughout
the rest of this section. Problem (\ref{eq:psse_ps_u}) remains nonconvex,
but we will show in Section \ref{sec:algorithm} that its specific
structure allows it to be approached as an eigenvector problem. While
that the reformulation from (\ref{eq:psse_subprob_x}) to (\ref{eq:psse_ps_u})
is not unique and can be done in a number of ways, the eigenvector
approach turns out to work exceptionally well for the \emph{positive
definite} and \emph{chordally sparse} $H$ that we propose in this
section.

The exactness of the reformulation from (\ref{eq:psse_subprob_x})
to (\ref{eq:psse_ps_u}) means that the optimal value $\opt(u)$ and
the optimal $x$ both coincide. This in turn requires the following
assumption, which effectively says that complex power is measured
in active-reactive pairs with the same standard deviation $\sigma_{i}$.

\begin{asm}\label{asm} Any complex bus and branch power measurements
come in real and reactive pairs with the same weight $w_{i}=w_{j}$.
\end{asm}

The basic rationale is to prevent scenarios where only one of the
two types of measurement is available, or where there is a significant
discrepancy in accuracy between the two. So long as the two precisions
$w_{i},w_{j}$ are on the order of same magnitude, the impact of approximating
them as exactly equal $w_{i}=w_{j}$ will be small. We emphasize that
in practical PSSE software packages, the weights $w_{i},w_{j}$ are
generally treated as hyperparameters, without strict adherence to
the actual statistical distribution of the associated measurements.

We begin by describing how an instance of the PSSE angle subproblem
(\ref{eq:psse_subprob_x}) with only complex bus power measurements
can be exactly reformulated into an instance of phase synchronization
(\ref{eq:psse_ps_u}). For simplicity, the proposition is stated to
cover all bus measurements, but any missing measurements can be easily
accommodated by setting the respective $w_{k}=0$. 
\begin{prop}[Bus complex power]
\label{prop:=0000201}Given $Y_{\bus}\in\C^{n\times n},$ $b\in\C^{n}$,
$u\in\R_{+}^{n},$ and $w\in\R_{+}^{n}$, the following holds 
\[
\sum_{k=1}^{n}w_{k}|v_{k}i_{k}^{*}-b_{k}|^{2}=x^{*}H_{\bus}(u)x\text{ for all }v=u\odot x
\]
where $i=Y_{\bus}v$ and $H_{\bus}(u)=C_{\bus}^{*}\diag(w)C_{\bus}$
and 
\[
C_{\bus}=\diag(u)Y_{\bus}\diag(u)-\diag(b^{*}).
\]
\end{prop}
\begin{IEEEproof}
Let $s=v\odot\conj(Y_{\bus}v)$ and $U=\diag(u)$ and note that 
\[
\conj(s)=\conj(v)\odot(Y_{\bus}v)=\diag(x^{*})UY_{\bus}Ux
\]
Define $\|z\|_{w}^{2}\equiv z^{*}\diag(w)z$ for $z\in\C^{n}$. Clearly,
\[
\sum_{k=1}^{n}w_{k}|v_{k}i_{k}^{*}-b_{k}|^{2}=\|s-b\|_{w}^{2}.
\]
Moreover, 
\begin{align*}
\|s-b\|_{w} & =\lVert\conj(s)-\conj(b)\|_{w}\\
 & =\lVert\diag(x^{*})UY_{\bus}Ux-\conj(b)\|_{w}\\
 & =\lVert\diag(x^{*})[UY_{\bus}Ux-\diag(x)\conj(b)]\|_{w}\\
 & =\lVert\diag(x^{*})C_{\bus}x\|_{w}=\lVert C_{\bus}x\|_{w}.
\end{align*}
The third line substitutes $\diag(x)\conj(b)=\diag(b^{*})x$. The
fourth line uses $\diag(x)\diag(x^{*})=I_{n}$. The final line is
because of the following 
\begin{align*}
\|\diag(x^{*})z\|_{w}^{2} & =z^{*}\diag(x)\diag(w)\diag(x^{*})z\\
 & =z^{*}\diag(w)z=\|z\|_{w}^{2}
\end{align*}
which uses the fact that diagonal matrices commute. 
\end{IEEEproof}
Next, we reformulate the PSSE angle subproblem (\ref{eq:psse_subprob_x})
with only complex branch power measurements into an instance of phase
synchronization (\ref{eq:psse_ps_u}). The proposition is stated with
``from'' measurements, though it also covers ``to'' measurements
with a trivial change of notations. Again, any missing measurements
can be accommodated by setting the respective $w_{k}=0$. 
\begin{prop}[Branch complex power]
\label{prop:=0000202}Given $Y_{\From}\in\C^{{l}\times n},$ $\alpha:[{l}]\to[n]$,
$b\in\C^{l}$, $u\in\R_{+}^{n},$ and $w\in\R_{+}^{{l}}$, the following
holds 
\[
\sum_{k=1}^{l}w_{k}|v_{\alpha(k)}i_{k,\From}^{*}-b_{k}|^{2}=x^{*}H_{\From}(u)x\text{ for all }v=u\odot x
\]
where $i_{\From}{=}Y_{\From}v$ and $D_{\From}=[\e_{k}\e_{\alpha(k)}^{T}]_{k=1}^{{l}}$
and $H_{\From}(u)=C_{\From}^{*}\diag(w)C_{\From}$ and 
\[
C_{\From}=\diag(D_{\From}u)Y_{\From}\diag(u)-\diag(b^{*})D_{\From}.
\]
\end{prop}
\begin{IEEEproof}
Let $s=(D_{\From}v)\odot\conj(Y_{\From}v)$ and $U=\diag(u)$ and
$U_{\From}=\diag(D_{\From}u)$ and $x_{\From}=D_{\From}x$ and note
that 
\[
\conj(s)=\conj(D_{\From}v)\odot(Y_{\From}v)=\diag(x_{\From}^{*})U_{\From}Y_{\From}Ux.
\]
Define $\|z\|_{w}^{2}\equiv z^{*}\diag(w)z$ for $z\in\C^{{l}}$.
Clearly, 
\[
\sum_{k=1}^{{l}}w_{k}|v_{\alpha(k)}i_{\From,k}^{*}-b_{k}|^{2}=\|s-b\|_{w}^{2}.
\]
Moreover, 
\begin{align*}
\|s-b\|_{w} & =\lVert\conj(s)-\conj(b)\|_{w}\\
 & =\lVert\diag(x_{\From}^{*})U_{\From}Y_{\From}Ux-\conj(b)\|_{w}\\
 & =\lVert\diag(x_{\From}^{*})[U_{\From}Y_{\From}Ux-\diag(x_{\From})\conj(b)]\|_{w}\\
 & =\lVert\diag(x_{\From}^{*})C_{\From}x\|_{w}=\lVert C_{\From}x\|_{w}.
\end{align*}
The third line substitutes $\diag(x_{\From})\conj(b)=\diag(b^{*})x_{\From}=\diag(b^{*})D_{\From}x$.
The fourth line uses $\diag(x_{\From})\diag(x_{\From}^{*})=I_{{l}}$.
The final line uses the fact that diagonal matrices commute. 
\end{IEEEproof}
Additionally, direct measurements of the voltage phasors $v_{i}$
can be expressed in this same form. Since phase is inherently relative,
all phasor measurements must be referenced to an underlying phasor
whose angle is assigned to zero by convention. Without loss of generality,
we assume that $v_{1}$ is the reference phasor, meaning that $\Im\{v_{1}\}=0$
by convention.

\begin{prop}[Phasor measurement]
Given $b\in\C^{n}$, $u\in\R_{+}^{n},$ and $w\in\R_{+}^{n}$, the
following holds 
\[
\sum_{i=1}^{n}w_{i}|v_{i}-b_{i}|^{2}=x^{*}H_{\pmu}(u)x\text{ for all }v=u\odot x,\Im\{v_{1}\}=0,
\]
where $H_{\pmu}(u)=(\diag(u)-\e_{1}b^{*})\diag(w)(\diag(u)-b\e_{1}^{T})$.
\end{prop}
\begin{IEEEproof}
Note that $\Im\{v_{1}\}=0$ implies that $x_{1}=1$. It follows that
$v_{i}-b_{i}=u_{i}x_{i}-b_{i}x_{1}=(u_{i}\e_{i}^{T}-b_{i}\e_{1}^{T})x$.
\end{IEEEproof}
Finally, let us consider the general PSSE angle subproblem (\ref{eq:psse_subprob_x})
with a mixture of voltage measurements, bus power measurements, branch
measurements, and phasor measurements. Evoking the propositions above
rewrites this problem as 
\begin{align*}
\min_{x\in\C^{n}}\quad & c(u)+x^{*}[H_{\bus}(u)+H_{\From}(u)+H_{\To}(u)+H_{\pmu}(u)]x\\
\text{s.t. }\quad & |x_{i}|=1\text{ for all }i\in\{1,\dots,n\}
\end{align*}
in which the contribution of the square-voltage measurements 
\[
c(u)=\sum_{i=1}^{n}w_{i,\mathrm{volt}}|u_{i}^{2}-b_{i,\mathrm{volt}}|^{2}
\]
is constant with respect to each $x_{i}=\exp(\i\theta_{i})=v_{i}/|v_{i}|$.
Therefore, up to a constant offset, a general instance of (\ref{eq:psse_subprob_x})
can be exactly reformulated into an instance of phase synchronization
(\ref{eq:psse_ps_u}) by setting $H(u)=H_{\bus}(u)+H_{\From}(u)+H_{\To}(u)+H_{\mathrm{pmu}}(u)$. 

In our experiments, we found that the efficacy of our approach is
critically driven by two spectral properties of the matrix $H(u)$:
the \emph{positive definiteness} of the matrix, and the size of the
\emph{spectral gap} between its two smallest eigenvalues. In turn,
these were driven by the \emph{observability} and the \emph{level
of noise} in the PSSE problem, as parameterized by $\lambda$ and
$\sigma$ defined in Theorem~1. Concretely, the matrix became \emph{positive
definite} $H(u)\succ0$ as soon as $\lambda>0$ and the PSSE problem
is observable, and the \emph{spectral gap} $\lambda_{n-1}(H(u))-\lambda_{n}(H(u))$
grew large whenever there was a large gap $\lambda-\sigma$ between
the observability and the level of noise. Based on these observations,
we speculate that our method begins to work as soon as it is possible
to solve the PSSE problem in an information-theoretic sense. But it
remains future work to provide a rigorous theoretical justification
for these empirical observations.

Due to space limitations, all of our numerical evaluation are performed
with a complete set of PQV bus measurements, i.e. we set $H(u)=H_{\bus}(u)$,
in order to ensure observability from a modest number of measurements.
The excellent results we obtained from this setup, in estimating
$2n-1$ unknowns from $m=3n$ noisy measurements, suggest that the
our approach does not require a majority of redundant measurements
to be effective. Nevertheless, another important future work is to
exhaustively study the impact of low observability or lack of observability
on the efficacy of our proposed approach.

\section{Solving Nonconvex Phase Synchronization}\label{sec:algorithm}

We now turn to the problem of solving the phase synchronization problem
(\ref{eq:psse_ps_u}), which we restate here without the dependence
on the voltage magnitude $u$: 
\begin{align}
\opt=\min_{x\in\C^{n}}\quad & x^{*}Hx\tag{P}\label{eq:P}\\
\text{s.t. }\quad & |x_{i}|=1\text{ for all }i\in\{1,\dots,n\}.\nonumber 
\end{align}
Recall that the $n\times n$ complex Hermitian matrix $H$, as previously
formulated in Section~\ref{sec:reformulation}, is generally positive
definite when the PSSE problem is observable. Another critical property
in our setting is that $H$ is \emph{chordally sparse}, meaning that
it can be symmetrically permuted and factored as $\Pi H\Pi^{T}=LL^{T}$
where the lower-triangular Cholesky factor $L$ is also sparse. Chordal
sparsity is prevalent across power applications, but is in fact much
stronger than the general sparsity seen in previous applications of
phase synchronization, like cryogenic electron microscopy~\cite{singer2011angular,boumal2016nonconvex}
and robot navigation~\cite{rosen2019se}. If the matrix $H$ is sparse
but not chordally sparse, then it can take sparse factorization up
to $O(n^{3})$ time and $O(n^{2})$ memory to solve the following
linear equations for some $y\in\R^{n}$ and $z\in\C^{n}$
\begin{equation}
[H-\diag(y)]v=z\text{ where }H-\diag(y)\succ0.\label{eq:mvp}
\end{equation}
These are the same worst-case complexity figures as if the matrix
$H$ were fully dense. But given that $H$ is chordal sparse in our
setting, the same equations (\ref{eq:mvp}) can be solved in guaranteed
$O(n)$ time and memory. This ability to solve (\ref{eq:mvp}) quickly
turns out to be an important application-specific structure that makes
solving phase synchronization substantially faster, more accurate,
and more robust.

In Section~\ref{subsec:init}, we show that the dominant eigenvector
of the matrix inverse $H^{-1}$ provides an initial guess $x^{(0)}$
known as a \emph{spectral initialization}, which we observe in our
experiments to be of exceptionally high quality. These initial angles
in $x^{(0)}$ can be iteratively refined, by parameterizing $x_{i}=\exp(\i\theta_{i})$
and applying a few steps of the Gauss--Newton method to the angles
$\theta$, to yield a globally optimal $x_{\opt}$ for (\ref{eq:P}),
which we recall is exactly equivalent to PSSE angle subproblem (\ref{eq:psse_subprob_x}).
Alternatively, multiplying $v^{(0)}=x^{(0)}\odot u$ by the magnitudes
$u$ used to calculate $x^{(0)}$ yields an initial guess for the
complete PSSE problem (\ref{eq:psse}), which can be iteratively refined
using classical Gauss--Newton. In Section~\ref{subsec:cert}, we
show how the minimum eigenvalue of a shifted matrix like $H-\diag(y)$
can be used to derive a formal mathematical proof, called a \emph{certification},
that a feasible choice of $x$ for (\ref{eq:P}) is $\delta$ near-globally
optimal. The estimated angles $\theta_{i}$ associated with a certifiably
globally optimal $x_{i}=\exp(\theta_{i})$ are guaranteed to be correct
(in the maximum likelihood sense), conditioned on the belief that
the estimated magnitudes $u$ are also correct.

While the basic ideas outlined above have been previously considered
for (\ref{eq:P}) arising in other applications, our main finding
in this paper is that they are \emph{unusually effective} for the
PSSE angle subproblem. We attribute this to three critical properties.
The most important is the \emph{chordal sparsity} of $H$, as it allows
us to use the inverse iteration to compute eigenvectors and eigenvalues,
instead of the (generalized) power iteration typically used for phase
synchronization. Second, the unusually high quality of the spectral
initialization $x^{(0)}$ in our setting is underpinned by the \emph{positive
definiteness} of our choice of $H$, which provides a strong baseline
guarantee on its quality. Third, we empirically observe that the instances
of (\ref{eq:P}) that arise from the PSSE angle subproblem have \emph{zero
duality gap}; it is this critical property that allows us to certify
global optimality exactly, without any conservatism.

\subsection{Spectral initialization}\label{subsec:init}

The spectral initialization for problem (\ref{eq:P}) is obtained
by computing the dominant eigenvector $v$ of the $n\times n$ positive
definite complex Hermitian matrix $H^{-1}$ 
\[
v=\arg\max_{v\ne0}\frac{v^{*}H^{-1}v}{v^{*}v}
\]
and then rescaling each element of this eigenvector onto the unit
circle 
\begin{align*}
x^{(0)}=\proj(v) & \equiv\arg\min_{x}\{\|v-x\|:|x|_{i}=1\text{ for all }i\}\\
 & =(v_{i}/|v_{i}|)_{i=1}^{n}.
\end{align*}
This choice of $x^{(0)}$ is derived by relaxing (\ref{eq:P}) into
the problem of computing the minimum eigenvalue of $H$ by summing
its $n$ constraints into a single constraint: 
\begin{align}
\opt & =\min_{x\in\C^{n}}\{x^{*}Hx:|x_{i}|^{2}=1\text{ for all }i\}\label{eq:=000020step1_1}\\
 & \ge\min_{v\in\C^{n}}\left\{ v^{*}Hv:v^{*}v=\sum_{i=1}^{n}|v_{i}|^{2}=n\right\} \label{eq:=000020step1_2}\\
 & =n\cdot\lambda_{\min}(H).\label{eq:=000020step1_3}
\end{align}
Despite the nonconvexity of the eigenvalue problem in (\ref{eq:=000020step1_2}),
it follows from the Courant--Fischer--Weyl Minimax Theorem that
the globally optimal $v$ is just the following eigenvector rescaled
to satisfy $\|v\|=\sqrt{n}$: 
\[
v=\arg\min_{v\ne0}\frac{v^{*}Hv}{v^{*}v}=\arg\max_{v\ne0}\frac{v^{*}H^{-1}v}{v^{*}v}.
\]
In general, this $v$ would not be feasible for (\ref{eq:P}), as
the fact that $\sum_{i=1}^{n}|v_{i}|^{2}=n$ does not necessarily
imply $|v_{i}|=1$ in reverse. However, it can be easily projected
onto the feasible set to yield a near-global optimal $x^{(0)}=\proj(v)$
for the original problem (\ref{eq:P}). Finally, we verify that $\proj(v)=\proj(\alpha v)$
for all $\alpha>0$, so it is unnecessary to rescale $v$ before projecting.

The fact that $H$ is positive definite provides a baseline global
optimality guarantee for $x^{(0)}$. Indeed, observe that (\ref{eq:=000020step1_1})-(\ref{eq:=000020step1_3})
implies the following lower-bounds on the cost of $x^{(0)}$: 
\begin{equation}
x^{(0)*}Hx^{(0)}\ge\opt\ge n\cdot\lambda_{\min}(H)>0.\label{eq:naive_bnd}
\end{equation}
If the voltage magnitudes are accurate and the noise level is small,
then we expect $\opt$ to be a small positive number. This, combined
with $\lambda_{\min}(H)>0$, ensures that the relaxation gap $\opt-n\cdot\lambda_{\min}(H)$
is also small. Therefore, the original problem (\ref{eq:=000020step1_1})
is close to being equivalent to its relaxation (\ref{eq:=000020step1_2}),
as in $\opt\approx n\cdot\lambda_{\min}(H)$, so we can expect projecting
the solution of the relaxation (\ref{eq:=000020step1_2}) back onto
to the original (\ref{eq:=000020step1_1}) to yield a near-global
solution $x^{(0)}$, with $x^{(0)*}Hx^{(0)}\approx n\cdot\lambda_{\min}(H)$.

For small-scale values of $n$ on the order of a few hundreds, the
dominant eigenvector $v$ of $H^{-1}$ can be reliably obtained by
computing the dense eigendecomposition $H=\sum_{i=1}^{n}\lambda_{i}v_{i}v_{i}^{*}$
for $\lambda_{1}\ge\cdots\ge\lambda_{n}$ in $\Theta(n^{3})$ time
and and $\Theta(n^{2})$ memory, and then taking $v\equiv v_{n}$.
For large-scale values of $n$ more than a thousand, it is better
to use the inverse iteration, which starts from an initial guess $z^{(0)}\in\C^{n}$
and iterates for $k=1,2,3\dots$ 
\[
z^{(k+1)}=(H-\alpha I)^{-1}z^{(k)}/\|(H-\alpha I)^{-1}z^{(k)}\|,
\]
where the shift factor $\alpha$ is chosen so that $H-\alpha I\succ0$,
and could be simply set to zero. As explained earlier, the chordal
sparsity of $H$ allows us to compute each $z^{(k)}\mapsto(H-\alpha I)^{-1}z^{(k)}$
in just $O(n)$ time and memory. In practice, we would precompute
the Cholesky factorization $\Pi(H-\alpha I)\Pi^{T}=LL^{T}$, and then
implement $(H-\alpha I)^{-1}z^{(k)}=\Pi^{T}L^{-T}L^{-1}\Pi z^{(k)}$
evaluated from the right to left. It follows from classical theory~\cite[Section~4.1]{saad2011numerical}
that the acute angle $\theta_{k}$ between the $k$-th iterate $z^{(k)}$
and $v$ converges as 
\[
\frac{\tan\theta_{k}}{\tan\theta_{0}}\le\left(\frac{\lambda_{n}-\alpha}{\lambda_{n-1}-\alpha}\right)^{k},\quad\tan\theta_{k}\equiv\frac{\|(I-vv^{*})z^{(k)}\|}{\|vv^{*}z^{(k)}\|},
\]
and that the error is $\|z^{(k)}-v\|=\sqrt{2-2\cos\theta_{k}}$. If
the real and imaginary parts of $z^{(0)}$ are sampled from the unit
Gaussian, then it takes at most $k=O(\frac{\lambda_{n}-\alpha}{\lambda_{n-1}-\lambda_{n}}\log(\frac{n}{\epsilon}))$
iterations to arrive at an $\epsilon$-accurate estimate of $v$ that
satisfies $\|z^{(k)}-v\|\le\epsilon$ with high probability.

\subsection{Global optimality certification}\label{subsec:cert}

Next, we construct a formal mathematical bound on the global optimality
of any feasible $x$ for problem (\ref{eq:P}). Its correctness hinges
on the following theorem. 
\begin{thm}[Lagrange duality]
\label{thm:dual}Let $\opt$ denote the optimal value of (\ref{eq:P}).
Then, 
\[
\opt\ge\max_{y\in\R^{n}}\;\one^{T}y+n\cdot\min\{0,\lambda_{\min}(H-\diag(y))\}.
\]
Moreover, if the above holds with equality, then the maximizer is
explicitly given as 
\[
y=\Re\{\conj(x_{\opt})\odot(Hx_{\opt})\},
\]
where $x_{\opt}$ denotes a global minimum to (\ref{eq:P}). 
\end{thm}
We defer the proof of Theorem~\ref{thm:dual} to the end of the subsection,
in order to explain how to use it to give a formally correct proof
of global optimality. Given $x\in\C^{n}$ satisfying $|x_{i}|=1$
for all $i$, we compute a \emph{lower-bound} $\mu$ on the following
eigenvalue problem 
\begin{equation}
\mu\le\lambda_{\min}[H-\diag(y)]\text{ where }y=\Re\{\conj(x)\odot(Hx)\}\label{eq:mu_def}
\end{equation}
that is ideally as tight (i.e. close to equality) as possible. It
then follows from $\one^{T}y=x^{*}Hx\ge\opt$ and Theorem~\ref{thm:dual}
that the global suboptimality of $x$ is bounded as: 
\begin{equation}
0\le x^{*}Hx-\opt\le\delta\quad\overset{\mathrm{def}}{=}\quad n\cdot\max\{0,-\mu\}.\label{eq:gap_def}
\end{equation}
We conclude that $x$ is \emph{certifiably} $\delta$ near-globally
optimal, because $x$ is feasible for the minimization problem (\ref{eq:P}),
and that its global suboptimality is guaranteed to be no worse than
$x^{*}Hx\le\opt+\delta$ by substituting $y$ into Theorem~\ref{thm:dual}.
This particular $y=\Re\{\conj(x)\odot(Hx)\}$ in turn serves as the\emph{
certificate} for the near global optimality claim, and the quantity
$\delta$ is called its \emph{duality gap}.

It is important to point out that zero-duality gap certificates with
$\delta=0$ do generally not exist. A globally optimal $x_{\opt}$
cannot usually be certified as \emph{exactly} globally optimal; the
best we can certify is that $x_{\opt}$ is \emph{nearly} globally
optimal for some strictly positive duality gap $\delta>0$ in (\ref{eq:gap_def}).
Indeed, this inherent conservatism reflects the general NP-hardness
of problem (\ref{eq:P}).

Therefore, it is surprising that zero-duality gap certificates seem
to always exist in the PSSE angle subproblem when the measurement
noise is reasonable. In our experiments, the spectral initialization
$x^{(0)}$ already yields a certificate $y^{(0)}$ with small duality
gap $\delta\approx0$. After parameterizing $x_{i}=\exp(\i\theta_{i})$
and performing just a single iteration of the Gauss--Newton method
on the angles $\theta$, the corresponding $y$ yields a duality gap
$\delta=0$ up to numerical roundoff. This empirical observation,
that global optimality in the PSSE angle subproblem can be exactly
certified with no conservatism, is the critical insight that underpins
the effectiveness of our approach.

In practice, a high-quality certification requires the lower-bound
$\mu$ to be as close to $\lambda_{\min}[H-\diag(y)]$ without exceeding
it. Unfortunately, any sparse eigenvector algorithm (e.g. the \texttt{eigs}
command in MATLAB/Octave) is only capable of computing an \emph{upper-bound}
like $\lambda_{\min}[H-\diag(y)]\le\frac{v^{*}[H-\diag(y)]v}{v^{*}v}$,
which will never give a correct certification. For small-scale values
of $n$, the lower-bound $\mu$ can be reliably obtained by computing
the dense eigendecomposition $H-\diag(y)=\sum_{i=1}^{n}\lambda_{i}v_{i}v_{i}^{*}$
for $\lambda_{1}\ge\cdots\ge\lambda_{n}$ in $\Theta(n^{3})$ time
and and $\Theta(n^{2})$ memory, and then shifting $\mu=\lambda_{n}-\epsilon$
by the numerical round-off $\epsilon\approx10^{-12}$ of the eigenvalue
algorithm. For large-scale values of $n$, we propose the inverse
bisection method in Algorithm~\ref{alg:bis}, which is based on the
fact that $H\succeq0$ if and only if it factorizes as $\Pi H\Pi^{T}=LL^{T}$.
The chordal sparsity of $H$ ensures that the cost of the Cholesky
factorization at each iteration costs $O(n)$ time, and termination
occurs in exactly $\left\lceil \log_{2}(\frac{\gamma-\mu}{\epsilon})\right\rceil $
iterations. The practical efficiency of the algorithm can be improved
by initializing the algorithm with a tighter bounds $\mu$ and $\gamma$,
e.g. by using a sparse eigenvector algorithm. In any case, the theoretical
complexity of an $\epsilon$-accurate certification is $O(n\log(1/\epsilon))$
time and $O(n)$ memory.

\begin{algorithm}
\caption{Inverse bisection method}\label{alg:bis}

\textbf{Input.} Chordally sparse $n\times n$ matrix $H\succ0$, fill-reducing
permutation $\Pi$ so that $\Pi H\Pi=LL^{T}$ factors into a sparse
Cholesky factor $L$, diagonal shift $y$, tolerance $\epsilon>0$.

\textbf{Output.} Lower-bound $\mu\le\lambda_{\min}(H-\diag(y))\le\mu+\epsilon.$

\textbf{Algorithm.} 
\begin{enumerate}
\item Initialize lower-bound $\mu$ and upper-bound $\gamma$ using any
$v\in\C^{n}$: 
\[
\mu=-\max_{i}\{0,y_{i}\},\quad\gamma=\frac{v^{*}[H-\diag(y)]v}{v^{*}v}.
\]
\item If $\gamma-\mu\le\epsilon$, return. Otherwise, set $\alpha=\frac{1}{2}(\mu+\gamma)$
and attempt to compute the Cholesky factor $L$ satisfying 
\[
\Pi[H-\diag(y)-\alpha I]\Pi^{T}=LL^{T}.
\]
\item If Cholesky factorization succeeds, set $\mu\gets\alpha$. Otherwise,
set $\gamma\gets\alpha$. Go back to Step 2. 
\end{enumerate}
\end{algorithm}

\begin{IEEEproof}[Proof of Theorem~\ref{thm:dual}]
The lower-bound follows by imposing a redundant constraint $\tr(xx^{*})\le n$,
relaxing the problem by substituting a semidefinite variable $X\succeq0$
in place of $xx^{*}$, and then taking the dual: 
\begin{align*}
\opt= & \min_{x\in\C^{n}}\{\inner H{xx^{*}}:\diag(xx^{*})=\one,\tr(xx^{*})\le n\},\\
\ge & \min_{X\succeq0}\{\inner HX:\diag(X)=\one,\tr(X)\le n\},\\
= & \max_{y\in\R^{n},\mu\in\R}\{\one^{T}y+n\mu:H-\diag(y)\succeq\mu I,\quad\mu\le0\}.
\end{align*}
The equality in the third line holds because Slater's condition is
satisfied. If equality is achieved in the second line, meaning that
$\inner HX=\one^{T}y+n\mu$ holds for a primal feasible $X=x_{\opt}x_{\opt}^{*}$
and dual feasible $S,y,\mu$ satisfying 
\begin{gather*}
\diag(X)=\one,\quad\tr(X)\le n,\quad X\succeq0,\\
\diag(y)+\mu I+S=H,\quad\mu\le0,\quad S\succeq0,
\end{gather*}
then $X,S,y,\mu$ must also satisfy the complementary conditions 
\begin{gather*}
SX=0,\quad\mu(\tr(X)-n)=0.
\end{gather*}
This immediately assures that $[H-\diag(y)]x_{\opt}=0$ holds. We
solve for $y$ in a least-squared sense (noting that $y$ needs to
be real) to yield 
\[
\min_{y\in\R^{n}}\|[H-\diag(y)]x_{\opt}\|=\Re\{\conj(x_{\opt})\odot(Hx_{\opt})\}.
\]
Note that 
\begin{align*}
 & \|[H-\diag(y)]x\|=\|Hx-\diag(x)y\|\\
= & \|\diag(x)(\diag(x^{*})Hx-y)\|=\|\conj(x)\odot(Hx)-y\|.
\end{align*}
To prove that $y$ uniquely satisfies $[H-\diag(y)]x_{\opt}=0$, suppose
that there exists some other $\tilde{y}\ne y$ such that $(H-\diag(y))x_{\opt}=(H-\diag(\tilde{y}))x_{\opt}=0$.
Then, subtracting the two equations yields $\diag(y-\tilde{y})x_{\opt}=\diag(x_{\opt})(y-\tilde{y})=0$,
and hence $y-\tilde{y}=0$ because $|x_{i,\opt}|=1$ and therefore
$\diag(x_{\opt})$ is a unitary matrix. 
\end{IEEEproof}

\section{Experimental Results}

Finally, we present experimental evidence to justify four empirical
claims made throughout the paper. In Section \ref{subsec:claim1},
we show that, where voltage magnitudes are accurate, the spectral
initialization is of \emph{exceptionally high quality}, to the extent
that further refinement is often unnecessary. In Section \ref{subsec:claim2},
we confirm that, starting with the spectral initialization, it usually
takes just a single Gauss--Newton iteration to arrive at certified
global optimality, with \emph{zero duality gap} and no conservatism.
In Section \ref{subsec:claim3}, we report that, as voltage magnitudes
become inaccurate, spectral initialization exhibits a \emph{graceful
degradation} in the accuracy of its estimated angles. Additionally,
the global optimality of the angles can still be certified, after
jointly refining the voltage magnitudes and angles with a few Gauss--Newton
iterations. Finally, in Section \ref{subsec:claim4}, we confirm that
to compute the initialization and certification takes \emph{time comparable
to between 2 and 6 Gauss--Newton iterations}.

Our first three experiments are performed to benchmark the accuracy
of estimating $v_{\gnd}\in\C^{n}$ from a complete set of bus power
$s_{k}$ and magnitude $|v_{i}|$ measurements, i.e. we estimate $2n-1$
unknowns from $m=3n$ measurements. In every experiment, the power
measurements are tainted with Gaussian noise $\epsilon_{i}\overset{\text{i.i.d.}}{\sim}\mathcal{N}(0,\sigma_{\text{noise}}^{2})$
with the same standard deviation $\sigma_{\text{noise}}$, and the
corresponding instance of PSSE (\ref{eq:psse}) is formulated with
unit weights $w_{1}=w_{2}=\cdots=w_{m}=1$. The first two experiments
consider perfect magnitude measurements, whereas the third experiment
allow the magnitudes to be noisy as well. The systems that we consider
cover the complete set of Polish, PEGASE, and RTE models from the
IEEE PES PGLib-OPF dataset~\cite{babaeinejadsarookolaee2019power}. 

The fourth experiment is performed to benchmark the worst-case running
time of our method. For this, we additionally include a complete set
of branch power $s_{k,\From}$ and $s_{k,\To}$ measurements. All
of our experiments are conducted on a workstation with 11th Gen Intel
Core i5-1135G7 2.40GHz processor and 8GB installed RAM. The code is
written in MATLAB (R2022a) and MATPOWER 7.1~\cite{zimmerman2010matpower}.

\subsection{Spectral initialization with accurate voltage magnitudes}\label{subsec:claim1}

This subsection demonstrates that spectral initialization provides
highly accurate initial voltage angle estimates for all Polish, PEGASE,
and RTE cases, for moderate noise in power measurements up to $\sigma_{\text{noise}}=0.1\,\text{pu}$,
assuming perfectly measured voltage magnitudes. We empirically establish
that the angular error of the initial guess is smaller than $178.3908\sigma_{\text{noise}}^{1.0013}$
degrees across all systems used, for $\sigma_{\text{noise}}\leq0.1\,\text{pu}$.
After one iteration of the Gauss-Newton algorithm, the angular error
is further reduced to $39.5507\sigma_{\text{noise}}^{1.0028}$ degrees,
for $\sigma_{\text{noise}}\leq0.1\,\text{pu}$. For example, with
$\sigma_{\text{noise}}=0.02\,\text{pu}$, the corresponding initial
angular error of $3.6^{\circ}$ is reduced to $0.8^{\circ}$ after
one Gauss-Newton iteration. In a less noisy scenario such as $\sigma_{\text{noise}}{=}0.005\,\text{pu}$,
the initial angular error starts at $0.9^{\circ}$ and refines to
$0.2^{\circ}$. However, given the already minimal initial error,
this improvement is practically negligible.
\begin{figure}[!t]
\includegraphics[width=0.48\textwidth]{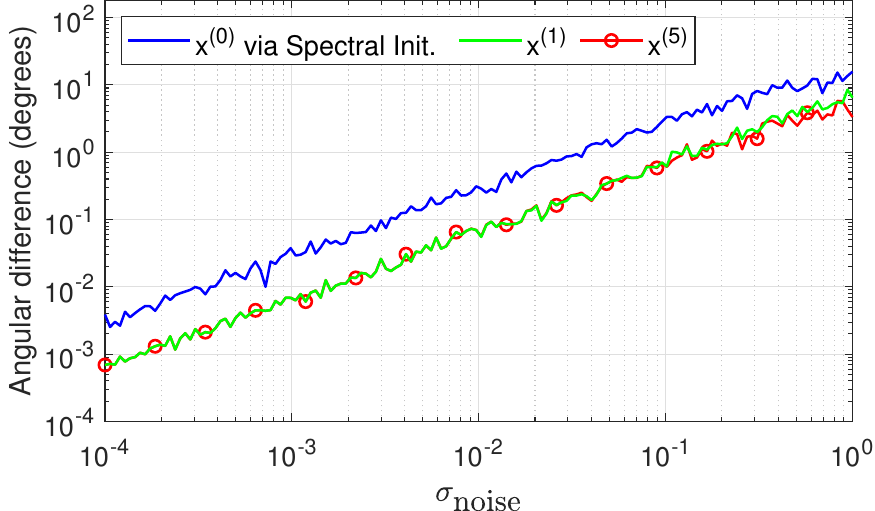}

\vspace{0.2cm}

\includegraphics[width=0.48\textwidth]{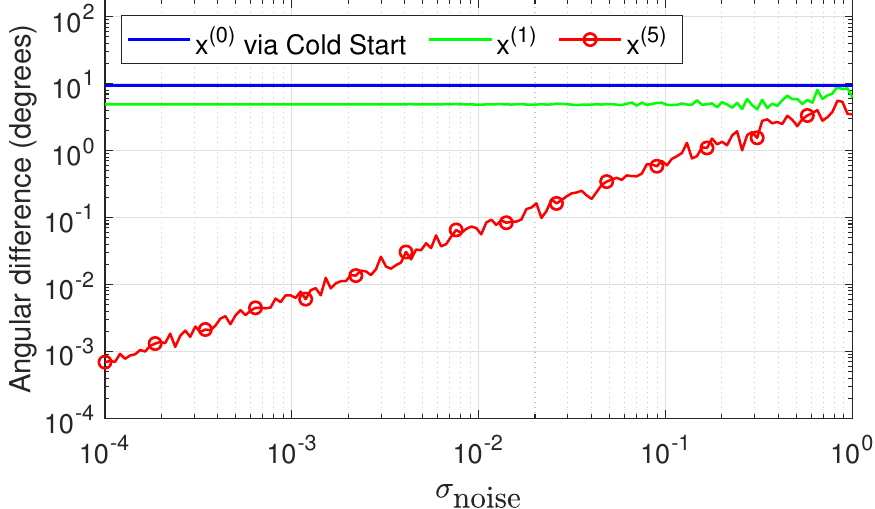}

\caption{\textbf{High quality of spectral initialization on 13k-bus system.}
(top) The initial guess chosen by spectral initialization is already
within $3^{\circ}$ of ground truth, for noise of up to $\sigma_{\mathrm{noise}}\le0.1$.
With just 1 iteration, GN converges to the MLE within $1^{\circ}$
of ground truth. (bottom) In contrast, GN with cold start takes 5
iterations to converge to the MLE.}

\label{fig:13k_angDiff}
\end{figure}

Beyond the first iteration, we observe that the angular error remains
essentially unchanged. This suggests that a single iteration of Gauss--Newton
algorithm is enough to locally optimizes the solution likelihood.
The consistent attainment of the same angular error level with a cold
start, even after multiple Gauss-Newton iterations, suggests that
the first iterate after spectral initialization had in fact been globally
optimal. In the next subsection, we will rigorously verify global
optimality up to $\sigma_{\text{noise}}=0.02\,\text{pu}$ using certified
lower bounds on the cost.

In our experiments, we solve the inner estimation problem over voltage
angles~(\ref{eq:psse_subprob_x}) with voltage magnitudes fixed as
the magnitude of the ground truth $u{=}|v_{\text{gnd}}|$ by setting
$\sigma_{\text{noise}}$ to a specific value. Then, we record the
maximum angular difference $\Delta_{\textrm{ang}}^{(i)}{=}\mathrm{max}\{|\theta_{k,\gnd}-\theta_{k}^{(i)}|,k=1,\dots,n\}$
in degrees, where $\theta_{\gnd}$ is the vector of ground truth angles
and $\theta_{k}^{(i)}$ for $x_{k}^{(i)}{=}\exp(\i\theta_{k}^{(i)})$
are the angles at the $k$-th bus associated with the $i$-th iterate.

First, we benchmark the initial guess for (\ref{eq:P}) using spectral
initialization and a cold start for angles, setting $\theta^{(0)}{=}0$
in Figure~\ref{fig:13k_angDiff}, which illustrates the average angular
difference $\Delta_{\textrm{ang}}$ across identical values of $\sigma_{\text{noise}}$
for the PEGASE 13659 bus system. In this figure, we conduct 50 random
trials per $\sigma_{\text{noise}}=10^{c}$ value, ranging from $10^{-4}$
to $10^{0}$ in 200 uniformly spaced increments over $c\in[-4,0]$.
The plots in Figure \ref{fig:13k_angDiff} (top) highlights a log-linear
relationship between the angular difference and the measurement noise,
expressed as $\Delta_{\textrm{ang}}^{(i)}{\approx}\beta\sigma_{\text{noise}}^{p},\,\beta,p\in\R$.

To establish upper bounds on angular errors, denoted by $\Delta_{\mathrm{ang}}^{(i)}\leq\beta^{(i)}\sigma_{\mathrm{noise}}^{p^{(i)}}$,
where $\beta^{(i)}$ and $p^{(i)}$ are positive parameters, we conduct
50 random trials for each $\sigma_{\text{noise}}$ value, ranging
from $10^{-3}$ to $10^{-1}$ in 150 uniformly spaced increments over
$c\in[-3,-1]$, across the Polish, PEGASE, and RTE case systems. Initially,
we identify the optimal $(\beta^{(i)},p^{(i)})$ parameters that best
fit the entire dataset, encompassing trials with $\sigma_{\text{noise}}$
in the range $[10^{-3},10^{-1}]$ for all systems. This involves employing
logarithmic regression analysis and solving a log-least squares problem~\cite{benoit2011linear},
yielding $(\beta^{(0)},p^{(0)}){=}(3.9946,1.0013)$ and $(\beta^{(1)},p^{(1)}){=}(1.7024,1.0028)$.
Subsequently, we adjust the leading terms $\beta^{(i)}$ to ensure
that the inequality $\Delta_{\mathrm{ang}}^{(i)}\leq\beta^{(i)}\sigma_{\mathrm{noise}}^{p^{(i)}}$
holds for all recorded samples, while retaining the optimal exponent
$p^{(i)}$. Following this adjustment, the leading terms $\beta^{(i)}$
are updated to $\beta^{(0)}=178.3908$ and $p^{(1)}=39.5507$. This
process establishes a valid upper bound on the observed angular error
across all case systems for $\sigma_{\text{noise}}$ in the range
$[10^{-4},10^{-1}]$, confirming the initial claim $\Delta_{\mathrm{ang}}^{(0)}\leq178.3908\sigma_{\text{noise}}^{1.0013}$
and $\Delta_{\mathrm{ang}}^{(1)}\leq39.5507\sigma_{\text{noise}}^{1.0028}$.

In addition, we give the median and maximum observed value of $\Delta_{\textrm{ang}}$
over $500$ random trials, setting $\sigma_{\text{noise}}=0.04$ pu,
in Table~\ref{tab:max_angle}.

\begin{table}[t]
\centering \caption{Median and max angular error in $500$ trials, $\sigma_{\text{noise}}=0.04$
pu}
\label{tab:max_angle} \vspace{-0.3cm}
\begin{minipage}[c]{0.99\columnwidth}%
 $\Delta_{\textrm{ang}}^{(i)}$ - maximum of angular error in degrees
for iteration $i$ starting from the spectral initialization. Numbers
in case names indicate the number of buses. %
\end{minipage} 

\medskip{}

\resizebox{\columnwidth}{!}{%
\begin{tabular}{|c|c|ccc|ccc|}
\hline 
\multirow{2}{*}{\#} & \multirow{2}{*}{Name} & \multicolumn{3}{c|}{Median $\Delta_{\textrm{ang}}$} & \multicolumn{3}{c|}{Max $\Delta_{\textrm{ang}}$}\tabularnewline
\cline{3-8}
 &  & $\Delta_{\textrm{ang}}^{(0)}$  & $\Delta_{\textrm{ang}}^{(1)}$  & \multicolumn{1}{c|}{$\Delta_{\textrm{ang}}^{(5)}$} & $\Delta_{\textrm{ang}}^{(0)}$  & $\Delta_{\textrm{ang}}^{(1)}$  & \multicolumn{1}{c|}{$\Delta_{\textrm{ang}}^{(5)}$}\tabularnewline
\hline 
1  & case1354pegase  & 0.13  & 0.08  & 0.08  & 0.6  & 0.39  & 0.39 \tabularnewline
2  & case1888rte  & 0.19  & 0.07  & 0.07  & 0.82  & 0.28  & 0.29 \tabularnewline
3  & case1951rte  & 0.18  & 0.08  & 0.08  & 0.73  & 0.37  & 0.37 \tabularnewline
4  & case2383wp  & 0.18  & 0.07  & 0.06  & 0.76  & 0.34  & 0.35 \tabularnewline
5  & case2736sp  & 0.15  & 0.08  & 0.08  & 0.68  & 0.38  & 0.38 \tabularnewline
6  & case2737sop  & 0.13  & 0.07  & 0.07  & 0.55  & 0.39  & 0.39 \tabularnewline
7  & case2746wop  & 0.13  & 0.07  & 0.07  & 0.59  & 0.28  & 0.28 \tabularnewline
8  & case2746wp  & 0.14  & 0.07  & 0.07  & 0.61  & 0.33  & 0.33 \tabularnewline
9  & case2848rte  & 0.15  & 0.09  & 0.09  & 0.68  & 0.37  & 0.38 \tabularnewline
10  & case2868rte  & 0.14  & 0.07  & 0.07  & 0.66  & 0.32  & 0.32 \tabularnewline
11  & case2869pegase  & 0.33  & 0.1  & 0.1  & 1.61  & 0.41  & 0.41 \tabularnewline
12  & case3012wp  & 0.16  & 0.08  & 0.08  & 0.78  & 0.37  & 0.37 \tabularnewline
13  & case3120sp  & 0.1  & 0.06  & 0.06  & 0.56  & 0.3  & 0.3 \tabularnewline
14  & case3375wp  & 0.16  & 0.09  & 0.09  & 0.72  & 0.41  & 0.41 \tabularnewline
15  & case6468rte  & 0.25  & 0.1  & 0.1  & 1.18  & 0.43  & 0.44 \tabularnewline
16  & case6470rte  & 0.22  & 0.09  & 0.09  & 1.13  & 0.45  & 0.45 \tabularnewline
17  & case6495rte  & 0.22  & 0.09  & 0.09  & 1  & 0.48  & 0.48 \tabularnewline
18  & case6515rte  & 0.25  & 0.08  & 0.08  & 1.05  & 0.4  & 0.4 \tabularnewline
19  & case9241pegase  & 1.33  & 0.17  & 0.17  & 6.56  & 0.8  & 0.81 \tabularnewline
20  & case13659pegase  & 1.18  & 0.26  & 0.26  & 5.32  & 1.33  & 1.33 \tabularnewline
\hline 
\end{tabular}} 
\end{table}

\subsection{Optimality certification with accurate voltage magnitudes}\label{subsec:claim2}

\begin{figure}[t!]
\centering \includegraphics[width=0.48\textwidth]{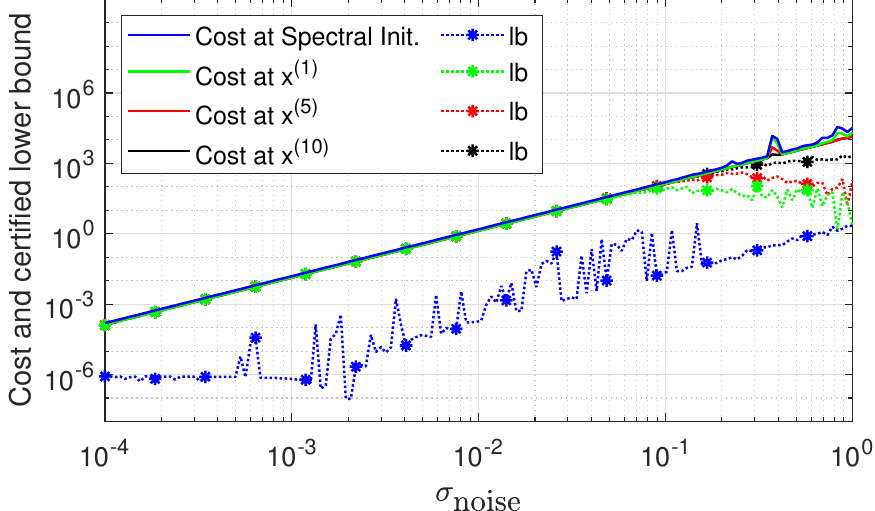}
\\
 \vspace{0.2cm}
 \includegraphics[width=0.48\textwidth]{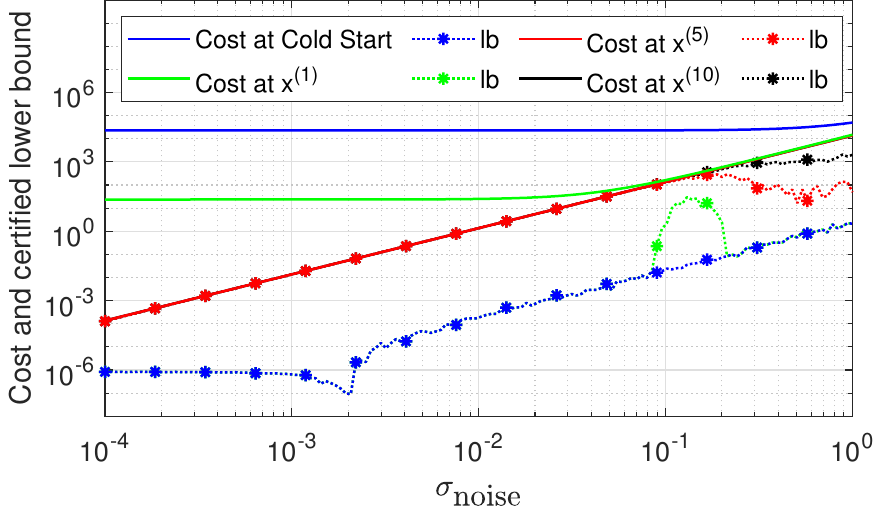}
\label{fig:13k_cost_lb_CS} \caption{\textbf{Certifiably optimal estimates on 13k-bus system using certified
lower bounds on the cost.} (top) It takes only one GN iteration starting
from the spectral initialization to certify global optimality on the
average case, for noise of up to $\sigma_{\mathrm{noise}}\le0.1$.
(bottom) The proposed method is applied to cold start initialization,
requiring 5 GN iterations to certify global optimality.}
\label{fig:13k_cost_lb} 
\end{figure}

This subsection presents numerical experiments demonstrating the efficacy
of certifying solutions under the assumption of perfectly observed
voltage magnitudes. The findings reveal that within a favorable noise
regime, up to $\sigma_{\text{noise}}=0.02\,\text{pu}$, we consistently
achieve a duality gap that is numerically zero to $10^{-6}$ relative
error, and certify global optimality to 99\% within just one iteration
from spectral initialization. Notably, across all Polish, PEGASE,
and RTE cases, within this noise range, solutions are certified to
100\% optimality within $10^{-5}$ relative error in 5 iterations
for any trial.

As $\sigma_{\text{noise}}$ increases, the duality gap widens, but
even with this increase, one iteration from spectral initialization
is typically sufficient to meet 99\% optimality certification. For
instance, with noise levels as high as $\sigma_{\text{noise}}=0.1\,\text{pu}$,
estimated angles after the first iteration still achieve a median
certification rate of 95\%. For this noise level, the subsequent iterations
still achieve a median rate of 99\% certification of global optimality.

The proposed certification method can be used with other initialization
methods. However, we observe that it takes considerably more iterations
to certify global optimality of a solution when a cold start initialization
is used since spectral initialization provides higher quality initial
estimates of the angles compared to a cold start approach. (See Figure~\ref{fig:13k_cost_lb}
for a comparison of the average cost and certified lower bound for
the PEGASE 13659 bus system.)

\begin{figure*}[!h]
\centering \includegraphics[width=0.95\textwidth]{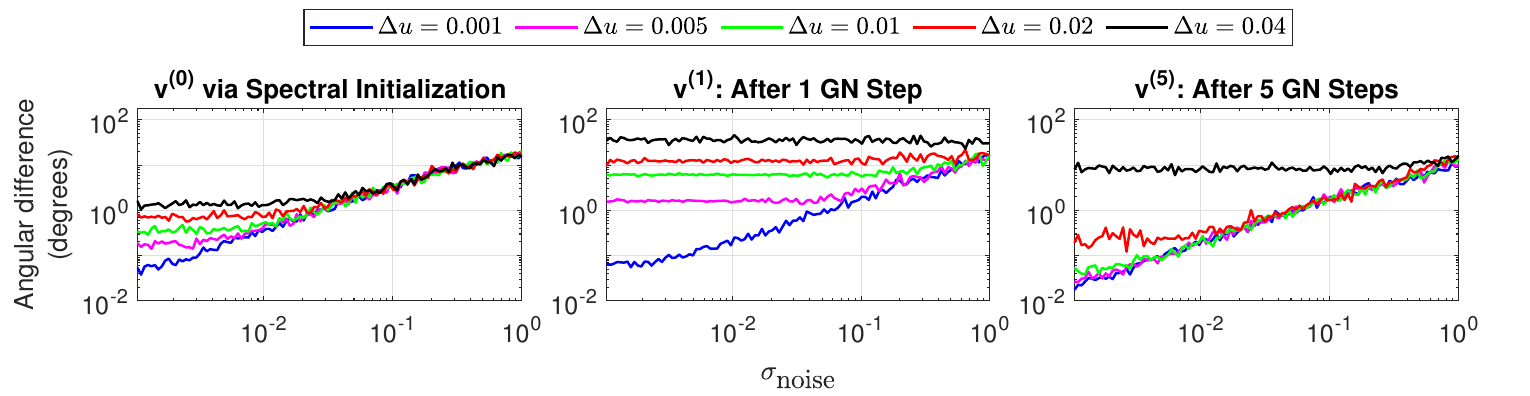}
\\
 \vspace{0.2cm}
 \includegraphics[width=0.95\textwidth]{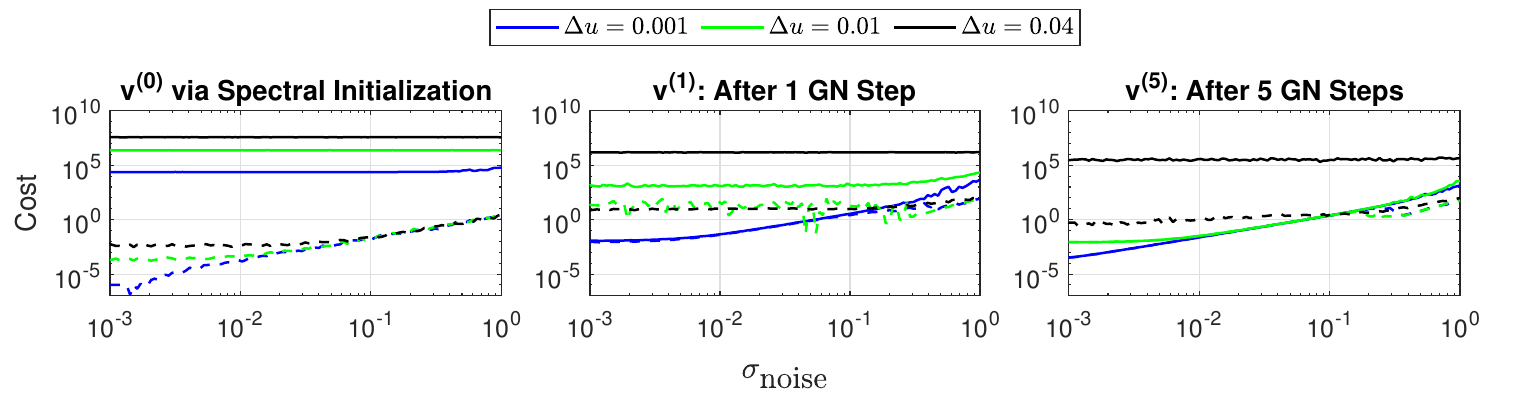}
\caption{\textbf{Initialization and optimality certification with inaccurate
voltage magnitudes on 13k-bus system.} (top) Spectral initialization
provides high quality initial guesses over the angles even with inaccurate
voltage measurements. In fact, spectral initialization yields more
accurate angle estimations compared to those achieved after 5 GN iterations
when $\Delta u=0.04$ pu. (bottom) With moderate voltage magnitude
error, the certification method verifies optimality of the angle estimates
once voltage magnitude estimates are accurate enough, after 5 GN iterations.}
\label{fig:Verr_13k} 
\end{figure*}

The simulation procedure in this subsection closely follows that of
the preceding subsection, focusing on evaluating the cost and the
certified lower bound for each trial. To reiterate, we define the
cost as $\mathrm{cost^{(i)}}=(x^{(i)})^{*}Hx^{(i)}$ and the certified
lower bound as $\mathrm{lb^{(i)}}\coloneqq\mathbf{1}^{T}y^{(i)}+n\cdot\min(0,\lambda_{\min}(H-\text{diag}(y^{(i)})))$,
where $y^{(i)}=\Re\{\conj(x^{(i)})\odot(Hx^{(i)})\}$. We present
the median and worst-case certified optimality over 500 random trials,
with $\sigma_{\text{noise}}=0.03$ pu, in Table~\ref{tab:min_opt}.
The values presented in the table are rounded to four digits after
the decimal point. Therefore, achieving $99.9999\%$ optimality in
this context is akin to achieving a duality gap that is effectively
zero within a relative error of $10^{-6}$.

\begin{table}[!t]
\centering \caption{Median and worst-case global optimality certification in $500$ random
trials, $\sigma_{\text{noise}}=0.03\,\text{pu}$}
\label{tab:min_opt} \vspace{-0.3cm}
 {%
\begin{minipage}[c]{0.99\columnwidth}%
 $\mathrm{cert}^{(i)}={\mathrm{cost}^{(i)}}/{\mathrm{lb}^{(i)}}$
: certified global optimality for iteration $i$ starting from the
spectral initialization, with values rounded to four digits after
the decimal.\\
\end{minipage}} \resizebox{\columnwidth}{!}{%
\begin{tabular}{|c|ccc|ccc|}
\hline 
\multicolumn{1}{|c|}{\multirow{2}{*}{\#}} & \multicolumn{3}{c|}{Median $\mathrm{cert}$ (\%)} & \multicolumn{3}{c|}{Min $\mathrm{cert}$ (\%)}\tabularnewline
\cline{2-7}
\multicolumn{1}{|c|}{} & $\mathrm{cert}^{(1)}$  & $\mathrm{cert}^{(5)}$  & \multicolumn{1}{c|}{$\mathrm{cert}^{(10)}$} & $\mathrm{cert}^{(1)}$  & $\mathrm{cert}^{(5)}$  & \multicolumn{1}{c|}{$\mathrm{cert}^{(10)}$}\tabularnewline
\hline 
1  & 99.9998  & 99.9999  & 99.9999  & 99.9971  & 99.9998  & 99.9999 \tabularnewline
2  & 99.9995  & 99.9999  & 99.9999  & 99.9894  & 99.9997  & 99.9997 \tabularnewline
3  & 99.9995  & 99.9999  & 99.9999  & 99.9934  & 99.9997  & 99.9997 \tabularnewline
4  & 99.9964  & 99.9999  & 99.9999  & 99.8590  & 99.9998  & 99.9997 \tabularnewline
5  & 99.9988  & 99.9999  & 99.9999  & 99.9839  & 99.9998  & 99.9998 \tabularnewline
6  & 99.9988  & 99.9999  & 99.9999  & 99.9807  & 99.9997  & 99.9997 \tabularnewline
7  & 99.9987  & 99.9999  & 99.9999  & 99.9677  & 99.9997  & 99.9996 \tabularnewline
8  & 99.9986  & 99.9999  & 99.9999  & 99.9653  & 99.9997  & 99.9997 \tabularnewline
9  & 99.9986  & 99.9999  & 99.9999  & 99.9914  & 99.9998  & 99.9998 \tabularnewline
10  & 99.9984  & 99.9999  & 99.9999  & 99.9830  & 99.9998  & 99.9998 \tabularnewline
11  & 99.9953  & 99.9999  & 99.9999  & 99.4636  & 99.9999  & 99.9999 \tabularnewline
12  & 99.9982  & 99.9999  & 99.9999  & 99.9647  & 99.9997  & 99.9997 \tabularnewline
13  & 99.9982  & 99.9999  & 99.9999  & 99.9392  & 99.9997  & 99.9997 \tabularnewline
14  & 99.9988  & 99.9999  & 99.9999  & 99.9867  & 99.9997  & 99.9997 \tabularnewline
15  & 99.9981  & 99.9999  & 99.9999  & 99.9684  & 99.9998  & 99.9998 \tabularnewline
16  & 99.9979  & 99.9999  & 99.9999  & 99.9776  & 99.9998  & 99.9998 \tabularnewline
17  & 99.9978  & 99.9999  & 99.9999  & 99.9898  & 99.9998  & 99.9998 \tabularnewline
18  & 99.9976  & 99.9999  & 99.9999  & 99.9850  & 99.9998  & 99.9998 \tabularnewline
19  & 99.9925  & 99.9999  & 99.9999  & 98.7865  & 99.0007  & 99.0007 \tabularnewline
20  & 99.9710  & 99.9999  & 99.9999  & 96.4800  & 96.6179  & 96.6179 \tabularnewline
\hline 
\end{tabular}} 
\end{table}

\subsection{Performance with inaccurate voltage magnitudes}\label{subsec:claim3}

Given complex power measurements, the phase synchronization formulation
proposed in Section~\ref{sec:reformulation} is exactly equivalent
to the PSSE problem under the assumption of perfectly measured voltage
magnitudes. Now, we demonstrate the performance of the spectral initialization,
and the ability to certify global optimality (subject to fixed voltages)
degrade gracefully with inaccurate voltage magnitudes. We empirically
establish that the angular error of initial estimates via spectral
initialization remains below $214.858\sigma_{\text{noise}}^{1.0938}+327.992\Delta u^{1.0153}$
degrees for $\sigma_{\text{noise}}\in[0.001,0.1]$ and $\Delta u\in[0.001,0.04]$,
where $\Delta u$ denotes the random voltage magnitude measurement
error introduced to the ground truth $u_{\gnd}$. For instance, with
$\sigma_{\text{noise}}=0.02$ and $\Delta u=0.02$, the angular error
is less than $9.2^{\circ}$. With more precise measurements (e.g.,
$\sigma_{\text{noise}}=0.005$ and $\Delta u=0.005$), the angular
error reduces below $2.2^{\circ}$.

Although the proposed method is particularly effective when initiated
with accurate voltage magnitudes, we observe that under conditions
of high voltage magnitude noise such as $\Delta u=0.04$ pu, spectral
initialization can yield more accurate angle estimations compared
to those achieved after numerous iterations of a local optimization
algorithm, as the inaccuracies in voltage magnitude measurements can
cause the algorithm to stagnate. This observation is illustrated in
Figure \ref{fig:Verr_13k} (top), which illustrates the average performance
across identical values of $\sigma_{\text{noise}}$ with fixed voltage
magnitude errors for the PEGASE 13659 bus system.

Moreover, the proposed certification method offers the capability
to validate the optimality of the angles of estimated complex bus
voltages $v^{(i)}$ at each iteration $i$, provided that the estimated
voltage magnitudes $u^{(i)}$ reach a satisfactory level of accuracy.
This aspect is illustrated in Figure \ref{fig:Verr_13k} (bottom).
The simulation procedure is similar to the preceding subsection, but
we introduce random error denoted as $\epsilon_{u}$ with the uniform
distribution $\epsilon_{u}\overset{\text{i.i.d.}}{\sim}\mathcal{U}([-\Delta u,\Delta u])$
to the actual values of voltage magnitudes. We apply the Gauss--Newton
algorithm to the PSSE problem in~(\ref{eq:psse}) to estimate both
bus voltage magnitudes and angles since the actual value of voltage
magnitudes are not known. Then, we record the angular error, cost,
and certified lower bound at each trial for the iterate $i$. The
angular error is defined as the maximum angular difference $\Delta_{\textrm{ang}}^{(i)}{=}\mathrm{max}\{|\theta_{k,\gnd}-\theta_{k}^{(i)}|,k=1,\dots,n\}$
in degrees, where $\theta_{\gnd}$ is the vector of ground truth angles
and $\theta_{k}^{(i)}$ for $v_{k}^{(i)}{=}u_{k}^{(i)}\exp(\i\theta_{k}^{(i)})$
are the angles at the $k$-th bus associated with the $i$-th iterate,
$u_{k}^{(i)}$ denoting the estimated voltage magnitude at the $k$-th
bus associated with the $i$-th iterate. We define the cost as $\mathrm{cost^{(i)}}=(x^{(i)})^{*}H(u^{(i)})x^{(i)}$
and the certified lower bound as $\mathrm{lb^{(i)}}\coloneqq\mathbf{1}^{T}y^{(i)}+n\cdot\min(0,\lambda_{\min}(H(u^{(i)})-\text{diag}(y^{(i)})))$,
where $y^{(i)}=\Re\{\conj(x^{(i)})\odot(H(u^{(i)})x^{(i)})\}$.

In Figure~\ref{fig:Verr_13k}, we conduct 50 random trials per $\sigma_{\text{noise}}=10^{c}$
value, ranging from $10^{-3}$ to $10^{0}$ in 150 uniformly spaced
increments over $c\in[-4,0]$, by setting $\Delta u$ to a specific
value. Initial guesses for angles are obtained via spectral initialization
for problem (\ref{eq:P}) with the voltage magnitudes fixed to $u^{(0)}=u_{\text{gnd}}+\epsilon_{u}$
and via cold start for angles with $\theta^{(0)}{=}0$ and $u^{(0)}=u_{\text{gnd}}+\epsilon_{u}$.
The certified lower bound and cost are calculated with the observed
$u^{(i)}$ after each iteration.

We undertake a procedure similar to the one outlined in Section~\ref{subsec:claim1}
to establish an upper bound on the initial angular error $\Delta_{\textrm{ang}}^{(0)}$.
To collect the necessary data, we introduce random noise to the power
measurements in 50 random trials per $\sigma_{\text{noise}}=10^{c}$
value, spanning from $10^{-3}$ to $10^{-1}$ in 150 uniformly spaced
increments across $c\in[-3,-1]$. In each set of 50 trials with a
specific $\sigma_{\text{noise}}$ value, we incorporate random error
into the actual voltage magnitude values with $\Delta u=10^{d}$,
varying from $0.001$ to $0.04$ in 50 uniformly spaced increments
over $d\in[-3,\log(0.04)]$. This results in a total of $150\times50$
trials for each of the Polish, RTE, and PEGASE case systems. We record
the angular error after the spectral initialization $\Delta_{\textrm{ang}}^{(0)}$
and associated $(\sigma_{\text{noise}},\Delta u)$ after each trial.

To model the relationship between measurement errors and angular error,
we employ a nonlinear function $\Delta_{\textrm{ang}}^{(0)}\approx\beta_{1}\sigma_{\text{noise}}^{p}+\beta_{2}\Delta{u}^{q}$.
After solving a log-least squares problem to determine the optimal
parameters $\boldsymbol{\beta}=(\beta_{1},\beta_{2},p,q)$ that best
fit the observations $\Delta_{\textrm{ang}}^{(0)}$ across all systems,
we find $\boldsymbol{\beta}=(4.4918,6.2656,1.0938,1.0153)$. Subsequently,
we adjust the leading terms $\beta_{1}$ and $\beta_{2}$ to ensure
that $\Delta_{\textrm{ang}}^{(0)}\leq\beta_{1}\sigma_{\text{noise}}^{p}+\beta_{2}\epsilon_{u}^{q}$
holds for all recorded samples, establishing the upper bound on initial
angle estimation error as $\Delta_{\textrm{ang}}^{(0)}\leq214.858\sigma_{\text{noise}}^{1.0938}+327.992\Delta u^{1.0153}$
degrees, as claimed.

\subsection{Cost of Computing Eigenvectors}\label{subsec:claim4}

\begin{table}[t]
\centering \caption{Timing for MATPOWER case studies when all bus and branch measurements
are used, $\sigma_{\text{noise}}=0.02$ p.u.}
\label{tab:times_all} \vspace{-0.15cm}
 {%
\begin{minipage}[c]{0.99\columnwidth}%
 \textbf{|$\mathcal{V}$| :} number of buses, \ m: total number of
measurements \ init: spectral initialization time~(millisecond),\ cert:
certification time~(millisecond), \ per-it: average time per GN
iteration~(millisecond) \vspace{0.15cm}
\end{minipage}} % \resizebox{\columnwidth}{!}{%
\begin{tabular}{|c|cc|ccc|cc|}
\hline 
\#  & $|\mathcal{V}|$  & m  & init  & cert  & per-it  & $\frac{\text{init}}{\text{per-it}}$  & $\frac{\text{cert}}{\text{per-it}}$ \tabularnewline
\hline 
1  & 1354  & 10672  & 17.9  & 9.7  & 6  & 3  & 1.6 \tabularnewline
2  & 1888  & 13900  & 24.6  & 13  & 7.8  & 3.2  & 1.7 \tabularnewline
3  & 1951  & 14286  & 28.2  & 15.1  & 8.8  & 3.2  & 1.7 \tabularnewline
4  & 2383  & 16350  & 62.1  & 34  & 16.6  & 3.7  & 2 \tabularnewline
5  & 2736  & 18548  & 71.7  & 38.8  & 19.7  & 3.6  & 2 \tabularnewline
6  & 2737  & 18550  & 70.6  & 38.8  & 20.5  & 3.4  & 1.9 \tabularnewline
7  & 2746  & 18720  & 73.3  & 37.2  & 20.1  & 3.6  & 1.9 \tabularnewline
8  & 2746  & 18608  & 66.9  & 34.9  & 18  & 3.7  & 1.9 \tabularnewline
9  & 2848  & 20800  & 67.2  & 35.2  & 21.2  & 3.2  & 1.7 \tabularnewline
10  & 2868  & 20968  & 67.3  & 37.5  & 20.7  & 3.3  & 1.8 \tabularnewline
11  & 2869  & 24066  & 72.1  & 39.9  & 22.9  & 3.1  & 1.7 \tabularnewline
12  & 3012  & 20312  & 75.9  & 40.2  & 20.5  & 3.7  & 2 \tabularnewline
13  & 3120  & 21012  & 80.7  & 40.2  & 20.5  & 3.9  & 2 \tabularnewline
14  & 3374  & 23392  & 88.4  & 46.3  & 22.9  & 3.9  & 2 \tabularnewline
15  & 6468  & 48936  & 210.8  & 93.8  & 45.7  & 4.6  & 2.1 \tabularnewline
16  & 6470  & 48960  & 205.4  & 93.4  & 44.4  & 4.6  & 2.1 \tabularnewline
17  & 6495  & 49066  & 208.2  & 92  & 45.4  & 4.6  & 2 \tabularnewline
18  & 6515  & 49178  & 211.3  & 94.7  & 46.4  & 4.6  & 2 \tabularnewline
19  & 9241  & 82678  & 388.9  & 174.4  & 82.9  & 4.7  & 2.1 \tabularnewline
20  & 13659  & 109186  & 541.4  & 245.4  & 107.2  & 5.1  & 2.3 \tabularnewline
\hline 
\end{tabular}% }
\end{table}

 Our proposed approach formulates the initialization and certification
as sparse eigenvector-eigenvalue problems. Our critical insight is
that the chordal sparsity of power systems allows us to solve this
eigenproblem efficiently using the inverse iteration, as explained
in Section~\ref{sec:algorithm}. Table~\ref{tab:times_all} demonstrates
that the average computation time for initialization is comparable
to 3-6 Gauss-Newton iterations, while the certification process takes
a time comparable to 2-3 Gauss-Newton iterations over 200 random trials
with $\sigma_{\text{noise}}=0.02$ pu, for all Polish, PEGASE, and
RTE cases encompassing all branch and bus power measurements. 

\section{Conclusions and future work}

This paper exactly reformulated the angle estimation subproblem of
PSSE as phase synchronization, thereby unlocking the powerful tools
of spectral initialization and global optimality certification for
PSSE. Our extensive experiments found that, for moderate levels of
noise, the spectral initialization is of such high quality that it
can sidestep the need for further iterative refinement in many practical
situations. Hence, spectral initialization has the potential to overcome
longstanding convergence issues in PSSE, including ill-conditioning
and spurious local minima. Moreover, with a single iteration of refinement,
the global optimality of the angles (with respect to fixed magnitudes)
are certified with zero duality gap, potentially providing an algorithmic
means of verifying correctness in PSSE. An important future work is
to provide a theoretical basis for our empirical observations described
above, with the goal of providing rigorous guarantees in line with
those given in~\cite{boumal2016nonconvex,bandeira2017tightness}. 

As another future direction, our method can be extended to work with
PQ measurements with different weights, hence overcoming the limitations
of Assumption~\ref{asm}. Consider the case in which the Q component
is completely unknown. We can introduce a pseudomeasurement over the
Q component (with the same weight as the P component), and then to
solve a secondary maximum likelihood estimation problem over the Q
component (with a different weight, reflecting its true level of unreliability).
If in a system with $n$ buses there exists only $p\ll n$ buses without
Q components, then this yields an inner-outer loop, in which the inner
subproblem over $n$ buses is easily solved using phase synchronization,
and the outer subproblem over the $p$ pseudomeasurements has far
fewer variables (and degree of nonconvexity). Since the angles are
most strongly impacted by the P component, the Q component does not
need to be accurately estimated before highly accurate estimations
of the angles can be obtained. 

\section*{Acknowledgments}

The authors thank Ian Hiskens, Amit Singer, and Andrew McRae for insightful
discussions and feedback, and two reviewers for comments that greatly
improved the paper. Financial support for this work was provided by
ONR Award N00014-24-1-2671, NSF CAREER Award ECCS-2047462, and by
C3.ai Inc. and the Microsoft Corporation via the C3.ai Digital Transformation
Institute.

{{\appendix[Three-Bus System]

It is well-known that the PSSE problem can exhibit spurious local
minima, that are similar the existence of low voltage solutions in
the classical power flow equations. But spurious local minima is a
\emph{strictly stronger} phenomenon~\cite{zhang2019spurious}, that
can exist even when the ambiguity between high- and low-voltage solutions
have been eliminated. Consider a fully-connected three-bus system,
with line impedances $z_{1,2}=\i0.03,z_{2,3}=\i0.08,z_{1,3}=\i0.03$
per unit, and zero charging susceptance. We take perfect voltage magnitude
measurements $u=(0.85,0.85,0.85)$~pu, and noisy real and reactive
power injection measurements $P=(-46.066,14.485,10.686)$~MW and
$Q=(273.22,238.86,153.4)$~MVa. The associated instance of PSSE seeks
to estimate the angles $\theta_{2},\theta_{3}$, with $\theta_{1}=0$
set as reference, and $u_{1},u_{2},u_{3}$ fixed. As shown in Figure~\ref{fig:3bus}
earlier, spurious local minima manifest even though the voltage magnitudes
are perfectly measured. }

\bibliographystyle{IEEEtran}
\bibliography{references}

\begin{IEEEbiography}[{\includegraphics[clip,width=1in,keepaspectratio,height=1.25in]{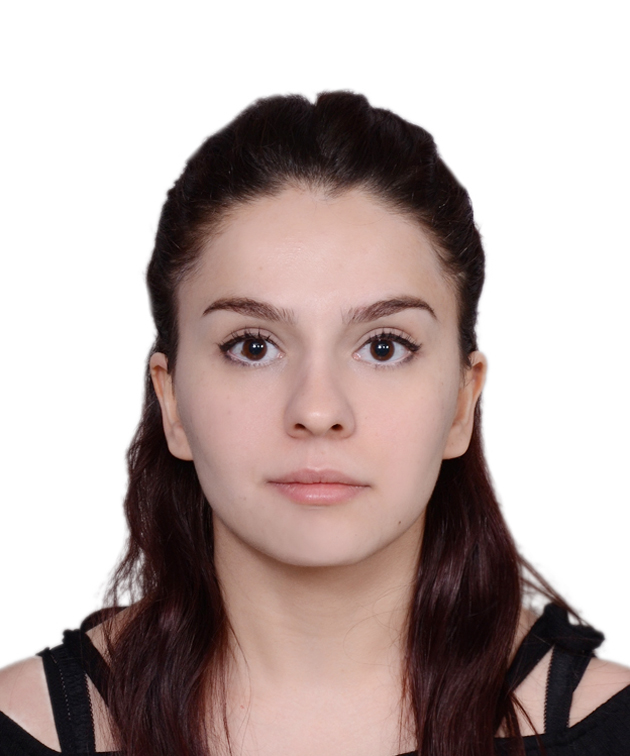}}]{Iven Guzel}
 (Student Member, IEEE) is a Ph.D. student at the Electrical and
Computer Engineering Department of the University of Illinois Urbana
Champaign (UIUC), USA. She received the B.S. and M.S. degrees in Electrical
and Electronics Engineering from Middle East Technical University
(METU), Turkey, in 2019 and 2022, respectively. 

Her research interests include practical algorithms for large-scale
systems and optimization with theoretical guarantees on cost, performance,
and reliability.
\end{IEEEbiography}

\begin{IEEEbiography}[{\includegraphics[width=1in]{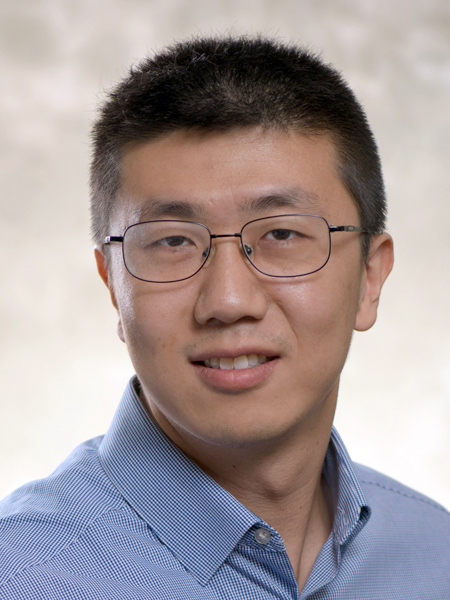}}]{Richard Y. Zhang}
 (Member, IEEE) received the B.E. (hons) degree with first class
honours in Electrical Engineering from the University of Canterbury,
Christchurch, New Zealand, in 2009 and the S.M. and Ph.D. degrees
in Electrical Engineering and Computer Science from the Massachusetts
Institute of Technology, Cambridge, MA, in 2012 and 2017 respectively.
From 2017 to 2019, he was a Postdoctoral Scholar in the Department
of Industrial Engineering and Operations Research at the University
of California, Berkeley, CA. 

He is currently an Assistant Professor in the Department of Electrical
and Computer Engineering at the University of Illinois, Urbana-Champaign.
His research interests are in mathematical optimization and its application
to power and energy systems. He received the NSF CAREER Award in 2021.
\end{IEEEbiography}

\end{document}